\documentclass[11pt]{amsart}
\usepackage{amsmath}
\usepackage{color}
\newtheorem{thm}{Theorem}[section]
\newtheorem{cor}[thm]{Corollary}
\newtheorem{lem}[thm]{Lemma}
\newtheorem{prop}[thm]{Proposition}
\theoremstyle{definition}

\theoremstyle{remark}
\newtheorem{rem}[thm]{\bf Remark}
\newtheorem{exe}[thm]{\bf Example}
\numberwithin{equation}{section}
\begin{document}
\title{Hypercyclic abelian semigroups of matrices on $\mathbb{R}^{n}$}

\author{Adlene Ayadi and Habib Marzougui}

\address{Habib Marzougui$^{1}$, University of Carthage, Faculty
of Science of Bizerte, (UR17ES21), ``Dynamical systems and their applications'',
7021, Jarzouna, Tunisia}
\email{habib.marzougui@fsb.rnu.tn}
 \address{Adlene Ayadi$^{2}$, \ University of Carthage, Faculty of Science of Bizerte, \\ 
(UR17ES21),``Dynamical systems and their
application'', University of Gafsa, Faculty of Science of Gafsa, Department of Mathematics, Tunisia}
 \email{adlenesoo@yahoo.com}

\thanks{This work is supported by the research unit: ``Dynamical systems and their applications'' (UR17ES21), Faculty of Science of Bizerte, Bizerte}


\begin{abstract}
 In this paper, we bring together results about the existence of a somewhere dense (resp. dense) orbit
 and the minimal number of generators for abelian semigroups of matrices on
$\mathbb{R}^{n}$. We solve the problem of determining the minimal number of matrices in normal form over 
$\mathbb{R}$ which form a hypercyclic abelian semigroup on 
$\mathbb{R}^{n}$. In particular, we show that no abelian semigroup generated by $\left[\frac{n+1}{2}\right]$
matrices on $\mathbb{R}^{n}$ can be hypercyclic. ($[\ ]$ denotes the integer part). 
This is a corrected version of the paper published in Topology and its Applications 210 (2016), 29--45 (see also  
\cite{aAhM16}). The differences between this version and the published version are explained
at the end of the Introduction.
%We provide a corrigendum to the statement and the proof of the implication $(ii) \Longrightarrow (iii)$ in the first part of Theorem %1.1., Lemma 3.6, Propositions 3.8 and 4.5.%and Corollary 1.8 
%and the proof of Theorems 1.1 and 1.4 of original article. Some others minors corrections are also shown. We also correct few typographical corrections. 
\end{abstract}
\maketitle

\section{\bf Introduction }
Let $M_{n}(\mathbb{R})$ be the set of all square matrices over
$\mathbb{R}$ of order $n\geq 1$ and by GL($n,\mathbb{R})$ the
group of invertible matrices of  $M_{n}(\mathbb{R})$. Let $G$ be
an abelian sub-semigroup of $M_{n}(\mathbb{R})$. For a vector
$v\in \mathbb{R}^{n}$, we consider the orbit of $G$ through $v$:
$G(v) = \{Av: \ A\in G\} \subset \mathbb{\mathbb{R}}^{n}$. The
orbit $G(v)\subset \mathbb{R}^{n}$ is called \textit{dense} (resp.
\textit{somewhere dense}) in ${\mathbb{R}}^{n}$ if $\overline{G(v)}=
{\mathbb{R}}^{n}$ (resp. $\mathring{\overline{G(v)}}\neq
\emptyset$),  where $\overline{E}$ (resp. $\overset{\circ}{E}$ )
denotes the closure (resp. the interior) of a subset $E\subset \mathbb{R}^{n}$. We say that $G$ is
 \textit{hypercyclic} if there exists a vector $v\in {\mathbb{R}}^{n}$ such that $G(v)$
is dense in ${\mathbb{R}}^{n}$. This concept was introduced by N. Feldman in \cite{fe1}. 
Recently there has been done much research on this subject. We mention in particular 
 [1, 2, 3, 4, 6, 7, 8, 12] for the abelian case and \cite{J} for the non-abelian case.
In the abelian case, Feldman \cite{fe1} showed that no semigroup generated by $n$-tuples of diagonalizable matrices on
$\mathbb{C}^{n}$ or $\mathbb{R}^{n}$ can be hypercyclic. 
If one removes the diagonalizability condition,
Costakis and al. \cite{chm2} proved that there is no semigroup generated by $n$-tuples
of non diagonalizable matrices on $\mathbb{R}^{n}$ which is
hypercyclic. Recently, Costakis and Parissis proved in
\cite{GC-IP} that the minimal number of matrices in Jordan form on
$\mathbb{R}^{n}$ which form a hypercyclic tuple is $n+1$. In the non-abelian case, Javaheri shows in \cite{J} that there exists a $2$-generator hypercylic semigroup in any
dimension in both real and complex cases. We refer the reader to
the recent book \cite{bm} and \cite{p} for a thorough account on
hypercyclicity. 

The aim of of this paper is to give necessary and sufficient conditions for $G$
to have a dense orbit (resp. a somewhere dense orbit).

First, we give a general result answering the above question for any abelian \textit{sub-semigroup} of $M_{n}(\mathbb{R})$ by 
providing an effective way of checking that a given semigroup is hypercyclic. Notice that in \cite{aAhM3}, the authors answer this question for any
abelian \textit{subgroup} of $\textrm{GL}(n;\mathbb{R})$,
so this paper can be viewed as a continuation of that work. We point out that the results obtained for groups 
are not sufficient for semigroups, the present paper is almost independent of \cite{aAhM3}.

 Second, we prove that the minimal number of matrices required to form a hypercyclic abelian semigroup having a 
normal form of length $(r+2s)$ is exactly 
$(n-s + 1)$ (see the definition below and Corollary \ref{C:5n}). In particular,
$\left[\frac{n+1}{2}\right]+ 1$ is the minimal number of matrices
on $\mathbb{R}^{n}$ required to form a hypercyclic abelian
semigroup on $\mathbb{R}^{n}$ ($[ \ ]$ denotes the integer part). This  answers a question raised by 
Feldman in (\cite{fe1}, Section 6). Notice that Abels and Manoussos \cite{hAaM}, Shkarin \cite{sh}
have, independently proved, similar results to Corollaries ~\ref{C:5n} and ~\ref{C:6n}.
The methods of proof in \cite{hAaM}, \cite{sh} and in this paper are quite different and have different 
consequences. 

To state our main results, we need to introduce the following
notations and definitions.

Set $\mathbb{N}$ be the set of non negative integers. Let $n\in\mathbb{N}$, $n\geq 1$ be fixed. For each
 $m = 1, 2, \dots, n,$ denote by:
\
\\
\textbullet \; $\mathbb{T}_{m}(\mathbb{R})$ the set of lower triangular matrices over
$\mathbb{R}$ with only one eigenvalue.
\
\\
\textbullet \; $\mathbb{T}_{m}^{*}(\mathbb{R}) =
\mathbb{T}_{m}(\mathbb{R})\cap \textrm{GL}(m, \mathbb{R})$ the
group of invertible matrices of $\mathbb{T}_{m}(\mathbb{R})$.
\
\\
\textbullet \; $\mathbb{T}_{m}^{+}(\mathbb{R})$ the group of
matrices of $\mathbb{T}_{m}(\mathbb{R})$ with all diagonal elements positive.
\
\\
\textbullet \; $\mathbb{S}$ the semigroup of matrices over $\mathbb{R}$ of
the form $\begin{bmatrix}
  \alpha & \beta \\
  -\beta & \alpha \\
  \end{bmatrix}.$
\
\\
\textbullet \; $\mathcal{B}_{0} = (e_{1},\dots,e_{n})$ the canonical basis of $\mathbb{R}^{n}$.
\
\\
\textbullet \; $I_{n}$ the identity matrix on $\mathbb{R}^{n}$.

 For each $1\leq m\leq \frac{n}{2}$, denote by
\\
\\
\textbullet \; $\mathbb{B}_{m}(\mathbb{R})$ the set of matrices of
$M_{2m}(\mathbb{R })$ of the form $$\begin{bmatrix}
  C &  &   & 0 \\
  C_{2,1} & C &   &   \\
  \vdots &  \ddots & \ddots & \\
  C_{m,1} & \dots & C_{m,m-1} & C
\end{bmatrix}: \  C, \ C_{i,j}\in \mathbb{S }, \ 2\leq i\leq m, 1\leq j\leq m-1.$$
\\
\textbullet \;  $\mathbb{B}^{*}_{m}(\mathbb{R}):=
\mathbb{B}_{m}(\mathbb{R})\cap \textrm{GL}(2m, \mathbb{R})$ the
group of invertible matrices of $\mathbb{B}_{m}(\mathbb{R})$.
\
\\
 
Let  $r, \ s\in \mathbb{N}$. By a partition of $n$ we mean a finite 
sequence of positive integers

$\eta = \begin{cases}
 (n_{1},\dots,n_{r};\  m_{1},\dots,m_{s}) & \textrm{ if } rs\neq  0, \\
(m_{1},\dots,m_{s}) & \textrm{ if } r= 0, \\ (n_{1},\dots,n_{r}) &
\textrm{ if } s=0
\end{cases}$
\
\\
such that $\underset{j=1}{\overset{r}{\sum}}n_{j} + 2\underset{j=1}{\overset{s}{\sum}}m_{j} =n$. In particular, we
have $r+2s\leq n$. The number 
$r+2s$ will be called the {\it length} of the partition. 
Given a partition $\eta =  (n_{1},\dots,n_{r};\  m_{1},\dots,m_{s})$, we denote by:
\
\\
\textbullet \; \ $\mathcal{K}_{\eta}(\mathbb{R}): =
\mathbb{T}_{n_{1}}(\mathbb{R})\oplus\dots \oplus
\mathbb{T}_{n_{r}}(\mathbb{R})\oplus \mathbb{B}_{m_{1}}(\mathbb{R
})\oplus\dots \oplus \mathbb{B}_{m_{s}}(\mathbb{R}).$

In particular:
\
\\
- $\mathcal{K}_{\eta}(\mathbb{R}) =
\mathbb{T}_{n}(\mathbb{R})$  and  $\eta=(n)$  if  $r=1, \ s=0$.
\
\\
- $\mathcal{K}_{\eta}(\mathbb{R}) =
\mathbb{B}_{m}(\mathbb{R})$  and  $\eta=(m)$,  $n=2m$ if $r=0, \ s=1$.
\
\\
- $\mathcal{K}_{\eta}(\mathbb{R}) =
\mathbb{B}_{m_{1}}(\mathbb{R})\oplus\dots\oplus\mathbb{B}_{m_{s}}(\mathbb{R})$
and $\eta=(m_{1},\dots,m_{s})$ if $r=0, \ s>1$.
\
\\
\textbullet \; $\mathcal{K}^{*}_{\eta}(\mathbb{R}): =
\mathcal{K}_{\eta}(\mathbb{R})\cap \textrm{GL}(n, \
\mathbb{R})$, it is a sub-semigroup of $\textrm{GL}(n, \
\mathbb{R})$.
\
\\
\textbullet \;  $\mathcal{K}^{+}_{\eta}(\mathbb{R}): =
\mathbb{T}^{+}_{n_{1}}(\mathbb{R})\oplus\dots \oplus
\mathbb{T}^{+}_{n_{r}}(\mathbb{R})\oplus
\mathbb{B}^{*}_{m_{1}}(\mathbb{R })\oplus\dots \oplus
\mathbb{B}^{*}_{m_{s}}(\mathbb{R}).$
\
\\

For a row vector $v\in\mathbb{R}^{n}$, we will be denoting by $v^{T}$ the transpose of $v$. We also have that
\
\\
\textbullet \;$u_{\eta} = [e_{\eta,1},\dots,e_{\eta,r}; f_{\eta,1},\dots,
f_{\eta,s}]^{T}\in \mathbb{R}^{n}$, where $e_{\eta,k} = [1,0,\dots,
0]^{T}\in \mathbb{R}^{n_{k}}$, $f_{\eta,l} = [1,0,\dots, 0]^{T}\in
\mathbb{R}^{2m_{l}}$, $ k=1,\dots, r; \ l=1,\dots, s$.
\
\\
\textbullet \;$f_{\eta}^{(l)} = [ 0,\dots,0,f^{(l)}_{1},\dots,
  f^{(l)}_{s}]^{T}\in \mathbb{R}^{n}$,
where 
\
\\
$f^{(l)}_{j} =
  \begin{cases}
    0\in \mathbb{R}^{2m_{j}} & \mathrm{if}\ \ j\neq l, \\
    [0,1,0,\dots,0]^{T}\in
\mathbb{R}^{2m_{l}} & \mathrm{if}\ \ j=l .
  \end{cases},  \ \ 1\leq l, j\leq s$
 \
 \\
Equivalently, $f_{\eta}^{(l)} = e_{t_{l}}$, where 
\
\\
$$t_{l} = \begin{cases}\underset{j=1}{\overset{r}{\sum}}n_{j}+2, \ & \textrm{ if } \ l=1\\
                 
     \underset{j=1}{\overset{r}{\sum}}n_{j}+
2\underset{j=1}{\overset{l-1}{\sum}}m_{j}+2, \ & \textrm{ if } \ ~l = 2,\dots,s\end{cases}$$

\
\\
Consider the matrix exponential map\; $\textrm{exp}: M_{n}(\mathbb{R})\longrightarrow \textrm{GL}(n, \mathbb{R})$ defined as $\textrm{exp}(M)=e^{M}$.
\
\\
Let $G$ be an abelian sub-semigroup of
$M_{n}(\mathbb{R})$. Then, following Proposition ~\ref{p:2}, there exists a $P\in \textrm{GL}(n, \mathbb{R})$ such that
$P^{-1}GP\subset \mathcal{K}_{\eta}(\mathbb{R})$ for some partition $\eta$ of $n$. Given two integers 
$r,s\in \mathbb{N}$ such that  $r+2s\leq n$, we shall say that the semigroup $G$ has ``a normal form of length 
$(r+2s)$'' if $G$ has a normal form 
in $\mathcal{K}_{\eta}(\mathbb{R})$ for some partition $\eta$ with length $(r+2s)$.
For such a choice of matrix $P$, we let:
\
\\
\textbullet \; $\mathrm{g}_{\eta}: = \textrm{exp}^{-1}(G)\cap \left[
P(\mathcal{K}_{\eta}(\mathbb{R}))P^{-1}\right]$.
\
\\
\textbullet \;$\mathrm{g}_{\eta}(u): = \{Bu: \ B\in \mathrm{g}_{\eta}\}, \
u\in\mathbb{R}^{n}$. 
\
\\
In particular, when $G\subset \mathcal{K}_{\eta}(\mathbb{R})$, then~
$\mathrm{g}_{\eta} = \textrm{exp}^{-1}(G)\cap
\mathcal{K}_{\eta}(\mathbb{R})$.
\
\\
\textbullet \;$G^{*}= G\cap \textrm{GL}(n,\mathbb{R})$.

\
\\
{\bf Definition} $($Index of $G$$)$. Set
$\widetilde{G^{*}} = P^{-1}G^{*}P$, where for every $M\in G^{*}$,
 $\widetilde{M} = P^{-1}M P\in \widetilde{G^{*}}$ which can be written as 
 $\widetilde{M} = \textrm{diag}(M_{1},\dots,M_{r}$; \ $\widetilde{M}_{1},\dots,\widetilde{M}_{s})\in
\mathcal{K}^{*}_{\eta}(\mathbb{R})$. Let $\mu_{k}$ be the eigenvalue
of $M_{k}$, $k=1,\dots,r$. We define the {\it index} of $\widetilde{G^{*}}$ to be
$$\textrm{ind}(\widetilde{G^{*}}):= \begin{cases}
0  \ \ \ \ \ \ \ \ \ \ \ \ \ \ \ \ \ \ \ \ \ \ \ \ \ \ \  \ \ \ \ \ \ \ \ \ \ \ \ \ \ \ \ \ \ \ \ \ \ \ \ \ \ \textrm{ if } \ r=0 \\
\\
\left\{\begin{array}{c} 1, \ \ \textrm{if} \ \ \ \mathrm{\exists}
\ \widetilde{M}\in \widetilde{G^{*}} \ \ \mathrm{with} \ \mu_{1}<0
\\ 0, \  \ \mathrm{otherwise}\ \ \ \ \ \ \ \ \ \ \ \ \ \ \ \ \ \ \
\ \ \
\end{array}
\right. \ \ \ \ \ \ \ \ \  \textrm{ if} \ \ r=1\\
\\
 \textrm{card}\left\{k\in\{1,\dots,r\}: \ \exists \widetilde{M}\in \widetilde{G^{*}} \ \textrm{ with } 
 \ \mu_{k} < 0, \ \mu_{i} >0, \ \forall\  i\neq k \right\} \ \textrm{if } \ r\notin\{0,\ 1\}.
 \end{cases}$$ 
\
\\

In particular,
\\
- If $\widetilde{G^{*}}\subset
\mathcal{K}^{+}_{\eta}(\mathbb{R})$ with $r\neq 0$, then
$\textrm{ind}(\widetilde{G})=0$.
\\
- If $\widetilde{G^{*}}\subset \mathbb{B}^{*}_{m}(\mathbb{R})$,
then $\textrm{ind}(\widetilde{G})=0$ (since $r=0$).\
\\

We define the {\it index} of $G$ to be $\textrm{ind}(G):
= \textrm{ind}(\widetilde{G^{*}})$. It is clear that this
definition does not depend on $P$.
\medskip

Our principal results can now be stated as follows:
\medskip

\begin{thm}\label{T:1n} Let $G$ be an abelian sub-semigroup of $M_{n}(\mathbb{R})$ and $P\in \textrm{GL}(n, \mathbb{R})$ such that $P^{-1}GP\subset
\mathcal{K}_{\eta}(\mathbb{R})$ where $\eta$ has length $(r+2s)$. 
\begin{enumerate}
 \item The following properties are equivalent:
 \begin{itemize}
  \item [(i)] $G$ has a somewhere dense orbit,
  \item [(ii)] The orbit $G(Pu_{\eta})$ is somewhere dense in $\mathbb{R}^{n}$,
  \item [(iii)] $\mathrm{g}_{\eta}(Pu_{\eta})$ is an additive sub-semigroup somewhere dense in $\mathbb{R}^{n}$.
\end{itemize}

\item Assume that $G$ is generated by $p$ matrices $A_{1},\dots,A_{p}$  ($p\geq 1$)
and let $B_{1},\dots,B_{p}\in \mathrm{g}_{\eta}$ such that $A_{1}^{2}=
e^{B_{1}},\dots, A_{p}^{2} = e^{B_{p}}$. Then $G$ has a somewhere dense orbit in $\mathbb{R}^{n}$ if and only if
$$\underset{k=1}{\overset{p}{\sum}}\mathbb{N}B_{k}Pu_{\eta} +
\underset{l=1}{\overset{s}{\sum}}2\pi\mathbb{Z}Pf_{\eta}^{(l)}$$ is dense in $\mathbb{R}^{n}$.
\end{enumerate}
\end{thm}

\begin{cor}\label{C:3} If $G$ is an abelian semigroup having a normal form of length $(r+2s)$ and 
generated by $(n-s)$ matrices of $M_{n}(\mathbb{R})$, then all orbits of $G$ are nowhere dense.
\end{cor}

\begin{cor}\label{C:4} If $G$ is an abelian semigroup generated by $[\frac{n+1}{2}]$ matrices of
$M_{n}(\mathbb{R})$, then all orbits of $G$ are nowhere dense.
\end{cor}
\medskip

\begin{cor} \label{c:1} Let $G$ be an abelian subgroup or a finitely generated semigroup 
of $M_{n}(\mathbb{R})$ and $P\in \textrm{GL}(n, \mathbb{R})$ such that $P^{-1}GP\subset
\mathcal{K}_{\eta}(\mathbb{R})$ where $\eta$ has length $(r+2s)$. 
\begin{enumerate}
 \item The following properties are equivalent:
 \begin{itemize}
  \item [(i)] $G$ has a somewhere dense orbit,
  \item [(ii)] The orbit $G(Pu_{\eta})$ is somewhere dense in $\mathbb{R}^{n}$,
  \item [(iii)] $\mathrm{g}_{\eta}(Pu_{\eta})$ is an additive sub-semigroup
dense in $\mathbb{R}^{n}$.
\end{itemize}

\item Assume that $G$ is an abelian subgroup generated by $p$ matrices $A_{1},\dots,A_{p}$  ($p\geq 1$)
and let $B_{1},\dots,B_{p}\in \mathrm{g}_{\eta}$ such that $A_{1}^{2}=
e^{B_{1}},\dots, A_{p}^{2} = e^{B_{p}}$. Then $G$ has a somewhere dense orbit in $\mathbb{R}^{n}$ if and only if
$$\underset{k=1}{\overset{p}{\sum}}\mathbb{Z}B_{k}Pu_{\eta} +
\underset{l=1}{\overset{s}{\sum}}2\pi\mathbb{Z}Pf_{\eta}^{(l)}$$ is dense in $\mathbb{R}^{n}$.
\end{enumerate}
\end{cor}

Corollary \ref{c:1} corresponds to Theorems 1.2 and 1.5 in \cite{aAhM3} when $G$ is an abelian subgroup.

\begin{rem} The implication $(ii)\Longrightarrow (iii)$ in the first part of Corollary \ref{c:1} is not true when $G$ is infinitely generated semigroup: Indeed, consider the semigroup $G = [-1,1]$ of $\mathbb{R}$. 
Here $G$ is infinitely generated, $u_{\eta} = 1$ and the orbit $G(u_{\eta}) = G$ is somewhere dense
in $\mathbb{R}$. However $\textrm{g}_{\eta}(u_{\eta})$ is not dense in $\mathbb{R}$ since $\textrm{g}_{\eta}(u_{\eta}) = \mathbb{R}_{-}$.
\end{rem}

\begin{thm}\label{T:2} 
Let $G$ be an sub-semigroup of $M_{n}(\mathbb{R})$ and $P\in \textrm{GL}(n, \mathbb{R})$ such that $P^{-1}GP\subset
\mathcal{K}_{\eta}(\mathbb{R})$ where $\eta$ has length $(r+2s)$.  
\begin{enumerate}
 \item The following properties are equivalent:
\begin{itemize}
  \item [(i)] $G$ is hypercyclic,
  \item[(ii)]  The orbit $G(Pu_{\eta})$ is dense in $\mathbb{R}^{n}$,
 \item [(iii)] $\mathrm{g}_{\eta}(Pu_{\eta})$ is an additive sub-semigroup dense in $\mathbb{R}^{n}$ and $\mathrm{ind}(G)=r$.
 \end{itemize}
 
 \item Assume that $G$ is generated by $p$ matrices $A_{1},\dots,A_{p}$  ($p\geq 1$)
and let $B_{1},\dots,B_{p}\in \mathrm{g}_{\eta}$ such that $A_{1}^{2}=
e^{B_{1}},\dots, A_{p}^{2} = e^{B_{p}}$. Then $G$ is hypercyclic if and only if 
 $\underset{k=1}{\overset{p}{\sum}}\mathbb{N}B_{k}Pu_{\eta} +
\underset{l=1}{\overset{s}{\sum}}2\pi\mathbb{Z}Pf_{\eta}^{(l)}$ is dense in $\mathbb{R}^{n}$ and
$\mathrm{ind}(G)=r$.
\end{enumerate}
\end{thm}
\medskip

 As a consequence of Theorem \ref{T:2}, we get in particular Theorem 1.3 and Corollary 1.6 of \cite{aAhM3}.
\medskip

\begin{thm}\label{T:4} For any partition $\eta$ of $n$ of length $(r+2s)$, there exist \\ $(n-s+1)$ matrices in
$\mathcal{K}_{\eta}^{*}(\mathbb{R})$ that generate a hypercyclic abelian semigroup.
\end{thm}
\medskip

From Theorem %~\ref{T:4n}, 
~\ref{T:4} and Corollary ~\ref{C:3}, % and \ref{C:4}, 
we obtain the following corollary.
\medskip

\begin{cor}\label{C:5n} The minimum number of matrices of $M_{n}(\mathbb{R})$
that generate an hypercyclic abelian semigroup with a normal form of length $(r+2s)$, is exactly $(n-s+1)$.
\end{cor}
\medskip

For $s=0$, we obtain:

\begin{cor}\label{C:8n}  The minimum number of trigonalizable matrices of $\mathbb{M}_{n}(\mathbb{R})$ that
generate a hypercyclic abelian semigroup is $n+1$.
\end{cor}
\medskip

In particular:

$-$ For $s=0$ and $r=n$, we obtain %Feldman's theorem
:
\medskip

\begin{cor} [%\cite{fe1}, Theorem 4.4; 
\cite{chm2}, Theorem 3.1] \label{C:7n}  The minimum number of diagonalizable matrices of $M_{n}(\mathbb{R})$ that generate a hypercyclic abelian semigroup is $n+1$.\
\end{cor}

$-$ For $s=0$ and $r = 1$, we obtain:

\begin{cor}\label{C:6nn} The minimum number of matrices of $\mathbb{T}_{n}(\mathbb{R})$ that generate an hypercyclic abelian semigroup, is $n+1$.
\end{cor}

\begin{cor}\label{C:6n} The minimum number of matrices of $M_{n}(\mathbb{R})$ that generate an hypercyclic abelian semigroup, is $\left[\frac{n+1}{2}\right]+1$.
\end{cor}
\medskip

This paper is organized as follows: In Section 2, we introduce the normal form of an abelian sub-semigroup of $M_{n}(\mathbb{R})$ and we give some related properties.
Sections 3 is devoted to the characterization of abelian
sub-semigroups of $\mathcal{K}^{*}_{\eta}(\mathbb{R})$ with a somewhere dense (resp. dense) orbit. In Section 4, we prove the first part of Theorems ~\ref{T:1n}, ~\ref{T:2} and Corollary \ref{c:1}.
Section 5 is devoted to finitely generated abelian semigroups; we prove the second part of
Theorems \ref{T:1n}, \ref{T:2}, and Corollary \ref{c:1}; Corollaries \ref{C:3} and \ref{C:4}. Theorem \ref{T:4} and
Corollary \ref{C:6n} are proved in Section 6. Section 7 is devoted to an
example for the case $n = 2$.
\medskip

\textit{Changes made to this version.}

We provide a corrigendum to the statement and the proof of the implication $(ii) \Longrightarrow (iii)$ in the first part of Theorem 1.1., Lemma 3.6, Propositions 3.8 and 4.5.%and Corollary 1.8 
and the proof of Theorems 1.1 and 1.4 of original article. Some others minors corrections are also down. We also correct few typographical corrections and polish some sentences. (See also \cite{aAhM16} for the corrigendum). 

The arguments in the proof of the implication $(ii) \Longrightarrow (iii)$ in the first part of Theorem 1.1 and the proof of Theorem 1.4 use Propositions 3.8 and 4.5 of original article. 
These later use Lemma 3.6 whose proof is not correct. So we give a correct statements and proofs of Lemma 3.6, Propositions 3.8 and 4.5. 
We also correct the proofs of the first part of Theorem 1.1 and the Theorem 1.4 while their statements are correct.
%Here we correct these results.

\section{\bf Normal form of abelian sub-semigroups of $M_{n}(\mathbb{R})$ and some related properties}
\medskip

First recall the following proposition.

\begin{prop}\label{p:1}$($\cite{aAhM3}, Proposition 1.1$)$ Let  $G$ be an abelian
\it{subgroup} of  $\textrm{GL}(n, \mathbb{R})$. Then there exists a
$P\in \textrm{GL}(n, \mathbb{R})$ such that $P^{-1}GP$  is an
abelian \it{subgroup} of $\mathcal{K}_{\eta}^{*}(\mathbb{R})$, for
some partition $\eta$ of $n$.
\end{prop}

The analogous proposition to Proposition ~\ref{p:1} for the sub-semigroup is the following.
\medskip

\begin{prop}\label{p:2}  Let  $G$ be an abelian
sub-semigroup of  $M_{n}(\mathbb{R})$. Then there exists a $P\in
GL(n, \mathbb{R})$ such that $P^{-1}GP$ is an abelian
sub-semigroup of $\mathcal{K}_{\eta}(\mathbb{R})$, for some partition $\eta$ of $n$.
\end{prop}

\begin{proof}  For every $A\in G$, there
exists $\lambda_{A}\in\mathbb{R}$ so that
$(A-\lambda_{A}I_{n})\in \textrm{GL}(n, \mathbb{R})$ (it suffices
to take $\lambda_{A}$ not an eigenvalue of $A$). Define
$\widehat{L}$ as the group generated by
$L:=\left\{A-\lambda_{A}I_{n}: \ A\in G\right\}$. Then
$\widehat{L}$ is an abelian subgroup of $\textrm{GL}(n,
\mathbb{R})$ and, by Proposition ~\ref{p:1}, there exists a $P\in
\textrm{GL}(n, \mathbb{R})$ such that $P^{-1}\widehat{L}P\subset
\mathcal{K}_{\eta}^{*}(\mathbb{R})$, for some
$\eta\in\mathbb{N}^{r+s}$ and $r,s\in \mathbb{N}$ with $r+2s\leq n$. As 
\
\\
$ P^{-1}LP = \left\{P^{-1}AP-\lambda_{A}I_{n}: \ A\in G\right\}$, we have
$P^{-1}GP\subset \mathcal{K}_{\eta}(\mathbb{R})$; this proves
the proposition.
\end{proof}
\medskip

The following results follow from basic properties of the matrix exponential map, and their proofs are left to the
reader.
\medskip

\begin{lem}\label{L:11} $\mathrm{exp}(\mathcal{K}_{\eta}(\mathbb{R})) =
\mathcal{K}^{+}_{\eta}(\mathbb{R})$.
\end{lem}

\begin{lem}\label{L:12} Let  $A, \ B\in \mathcal{K}_{\eta}(\mathbb{R})$. If  $e^{A}e^{B} =
e^{B}e^{A}$, then $AB = BA$.
\end{lem}

\begin{lem}\label{L:05} $($\cite{wR}, Proposition $7^{\prime}$, page 17$)$. The restriction
\
\\
$\mathrm{exp}_{|\mathbb{T}_{n}(\mathbb{R})}: \mathbb{T}_{n}(\mathbb{R})\longrightarrow\
\mathbb{T}_{n}^{*}(\mathbb{R})$ is a local diffeomorphism, in particular
it is an open map.
\end{lem}
\smallskip

\begin{cor}\label{LLL:1} The restriction $\mathrm{exp}_{|\mathcal{K}_{\eta}(\mathbb{R})}: \mathcal{K}_{\eta}(\mathbb{R})\longrightarrow
\mathcal{K}^{*}_{\eta}(\mathbb{R})$ is a local diffeomorphism, in particular it is an open map.
\end{cor}
\smallskip

\begin{proof} The proof results from Lemma ~\ref{L:05} and the fact that
\
\\
$\textrm{exp}_{|\mathcal{K}_{\eta}(\mathbb{R})} 
= \textrm{exp}_{|\mathbb{T}_{n_{1}}(\mathbb{R})}\oplus\dots\oplus \textrm{exp}_{|\mathbb{T}_{n_{r}}(\mathbb{R})}\oplus \textrm{exp}_{|\mathbb{B}_{m_{1}}(\mathbb{R})}\oplus\dots\oplus \textrm{exp}_{|\mathbb{B}_{m_{s}}(\mathbb{R})}$.
\end{proof}
\medskip

We let \; $$U : = \underset{k=1}{\overset{r}{\prod}}(\mathbb{R}^{*}\times\mathbb{R}^{n_{k}-1})\times
\underset{l=1}{\overset{s}{\prod}}\left((\mathbb{R}^{2}\backslash\{(0,0)\})\times\mathbb{R}^{2m_{l}-2}\right)$$

and \; $$C_{u_{\eta}} = \underset{k=1}{\overset{r}{\prod}}(\mathbb{R}^{*}_{+}\times\mathbb{R}^{n_{k}-1})\times
\underset{l=1}{\overset{s}{\prod}}\left((\mathbb{R}^{2}\backslash\{(0,0)\})\times\mathbb{R}^{2m_{l}-2}\right)$$

\
\\
where $u_{\eta}$ is defined in the introduction back on page $3$.
It is plain that $U$ is open and dense in $\mathbb{R}^{n}$ and that $C_{u_{\eta}}$ is the connected component of $U$  containing 
$u_{\eta}$.

\begin{lem}\label{L:74} $($\cite{aAhM3}, Lemma 3.3$)$ Let $u\in U$.
\
\\
 $($i$)$ If $B\in \mathcal{K}_{\eta}(\mathbb{R})$ satisfies $Bu\in U$, then
$B\in \mathcal{K}^{*}_{\eta}(\mathbb{R})$.\
\\
$($ii$)$ If $B\in \mathcal{K}_{\eta}(\mathbb{R})$ satisfies $Bu\in  C_{u_{\eta}}$, then
$B\in \mathcal{K}^{+}_{\eta}(\mathbb{R})$.
\end{lem}

\section{\bf Abelian sub-semigroup of
$\mathcal{K}^{*}_{\eta}(\mathbb{R})$ with a somewhere dense (resp. dense)
orbit}

Throughout this section, we let $G$ be an abelian sub-semigroup of
$\mathcal{K}^{*}_{\eta}(\mathbb{R})$. Define
\
\\
\textbullet \; $\mathcal{C}(G):= \{A\in
\mathcal{K}_{\eta}(\mathbb{R}): \ AB=BA, \ \forall \ B\in G
\}.$

Since $G$ is abelian, $G\subset \mathcal{C}(G)$.
\
\\
\textbullet \; $G^{+}: = G\cap \mathcal{K}^{+}_{\eta}(\mathbb{R})$.
\medskip

\begin{lem}\label{L:1} We have that\
\begin{itemize}
%\item[(i)]  $\mathcal{C}(\widehat{G}) = \mathcal{C}(G)$.
\item[(i)] $\mathrm{g}_{\eta} \subset \mathcal{C}(G)$ and all matrices of $\mathrm{g}_{\eta}$ commute, 
\item[(ii)] $\mathrm{exp}(\mathrm{g}_{\eta})=G^{+}$,
\item[(iii)] $\mathrm{exp}(\mathcal{C}(G))=\mathcal{C}(G)\cap
\mathcal{K}^{+}_{\eta}(\mathbb{R})$.
\end{itemize}
\end{lem}

\begin{proof} %(i) If $B\in \mathcal{C}(G)$ and $A\in G$ then $A^{-1}B=BA^{-1}$ (since $AB=BA$). We conclude that $B\in\mathcal{C}(\widehat{G})$.\
%\\
(i) By Lemma ~\ref{L:12}, all elements of  $\mathrm{g}_{\eta}$  commute, hence
 $\textrm{g}_{\eta}\subset \mathcal{C}(\textrm{g}_{\eta})$. Let $B\in\mathrm{g}_{\eta}$ and $A\in G$. 
 So $e^{B}\in G^{+}\subset G$. As $G$ is abelian,
 $Ae^{B}=e^{B}A$. Hence $e^{A}e^{B}=e^{B}e^{A}$. Since $A, \ B\in \mathcal{K}_{n}(\mathbb{R})$, it follows
 by Lemma ~\ref{L:12} that  $AB = BA$ and therefore $B\in \mathcal{C}(G)$. We conclude that $\mathrm{g}_{\eta}\subset
\mathcal{C}(G)$.
\\
(ii) We have $\textrm{exp}(\textrm{g}_{\eta})\subset G^{+}$ by definition. Conversely, let  $A\in G^{+}$. 
By Lemma \ref{L:11}, there exists
$B\in\mathcal{K}_{\eta}(\mathbb{R})$  such that  $e^{B} = A$. Hence $B\in \textrm{exp}^{-1}(G)\cap
\mathcal{K}_{\eta}(\mathbb{R}) = \mathrm{g}_{\eta}$, and then  $A\in
\textrm{exp}(\mathrm{g}_{\eta})$. So $G^{+}\subset
\textrm{exp}(\mathrm{g}_{\eta})$, this proves (ii).
\
\\
(iii) Let $A=e^{B}$, where $B\in \mathcal{C}(G)$, and let $C\in
\mathcal{C}(G)$. Then $BC=CB$, and therefore $Ce^{B}=e^{B}C$, or $AC=CA$. It follows that $A\in \mathcal{C}(G)$. Since $B\in
\mathcal{K}_{\eta}(\mathbb{R})$, so  $A\in
\mathcal{K}_{\eta,r,s}^{+}(\mathbb{R})$, by Lemma~\ref{L:11}.
Conversely, let $A\in \mathcal{C}(G)\cap
\mathcal{K}^{+}_{\eta}(\mathbb{R})$. By Lemma~\ref{L:11},
there exists $B \in\mathcal{K}_{\eta}(\mathbb{R})$ such that
$e^{B} = A$. Let $C\in G$. Then $Ce^{B} = e^{B} C$, and hence
$e^{C}e^{B} = e^{B} e^{C}$. Since $B, C
\in\mathcal{K}_{\eta}(\mathbb{R})$, it follows by
Lemma~\ref{L:12} that $BC = CB$. Therefore, $B \in \mathcal{C}(G)$,
and hence $A \in \textrm{exp}(\mathcal{C}(G))$.
\end{proof}
\smallskip

Let $u \in \mathbb{R}^{n}$ and consider the linear map 
\begin{align*}\Phi_{u} :
\mathcal{C}(G)\longrightarrow \mathbb{R}^{n}\;\;A\longmapsto Au
\end{align*}

Denote by Vect$(G)$ the vector subspace of $\mathcal{K}_{\eta}(\mathbb{R})$ generated by $G$.
\
\\
\begin{prop}\label{pp:01} If\; $\overset{\circ}{\overline{G(u)}}\neq\emptyset$ (resp. 
$\overset{\circ}{\overline{\mathrm{g}_{\eta}(u)}}\neq\emptyset$) for some $u
\in \mathbb{R}^{n}$, then $\Phi_{u}$ is a linear isomorphism.
Moreover, $\Phi_{u}(\mathrm{Vect}(G))= \Phi_{u}(\mathcal{C}(G))=\mathbb{R}^{n}$.
\end{prop}

\begin{proof} Case 1: $\overset{\circ}{\overline{G(u)}}\neq\emptyset$.
\
\\
 - $\Phi_{u}$ is surjective: we have that
$\Phi_{u}(\mathcal{C}(G))$ is a vector subspace of
$\mathbb{R}^{n}$. Since $G\subset \mathcal{C}G)$, it follows that
$G(u)\subset \Phi_{u}(\mathcal{C}(G))$. As $\mathcal{C}(G)$ is a
vector space, $\emptyset\neq \overset{\circ}{\overline{G(u)}}\subset
 \Phi_{u}(\mathcal{C}(G))$, and therefore $\Phi_{u}(\mathcal{C}(G))=\mathbb{R}^{n}$.
  We also have $G(u)\subset \Phi_{u}(\textrm{Vect}(G))$, so as above, 
  $\Phi_{u}(\textrm{Vect}(G)) = \mathbb{R}^{n}$.
\
\\ 
- $\Phi_{u}$ is injective: let $A \in \textrm{Ker}(\Phi_{u})$, so $Au=0$.
  Let $x\in \mathbb{R}^{n}$; then, by the above working there exists
  $B\in \textrm{Vect}(G)$ such that $x=Bu$. As $A\in \textrm{Ker}(\Phi_{u})\subset
  \mathcal{C}(G)$, then $AB= BA$. Therefore, $Ax=ABu= BAu=
  B(0)=0$. It follows that $A= 0$, and hence
  $\textrm{Ker}(\Phi_{u})=\{0\}$.
\
\\
 Case 2: \; $\overset{\circ}{\overline{\mathrm{g}_{\eta}(u)}}\neq\emptyset$. We 
 also have that $\Phi_{u}(\mathcal{C}(G))=\mathbb{R}^{n}$ since $\mathrm{g}_{\eta}\subset \mathcal{C}(G)$ (Lemma
~\ref{L:1}, \ (i)) and $\emptyset\neq \overset{\circ}{\overline{\mathrm{g}_{\eta}(u)}}\subset \Phi_{u}(\mathcal{C}(G))$.
\end{proof}

A subset $E\subset \mathbb{R}^{n}$ is called \emph{$G$-invariant}
if $A(E)\subset E$  for any  $A\in G$. In particular, if $E$
is $G$-invariant, then so is $\overline{E}$ and $\mathring{E}$.
\medskip

\begin{lem}\label{L:00} We have that
\begin{itemize}
\item [(i)] $\mathbb{R}^{n}\backslash U$ is a union of $(r+s)$ $G$-invariant vector subspaces of
$\mathbb{R}^{n}$ of dimension $n-1$ or $n-2$.
\item [(ii)] $U$ and $\mathbb{R}^{n}\backslash U$ are $G$-invariant.
\end{itemize}
 \end{lem}

\begin{proof} Assertion (i) is trivial. Assertion (ii): A simple calculation from the definition yields
that $U$ and $\mathbb{R}^{n}\backslash U$ are $G$-invariant.
 \end{proof}

\begin{lem}\label{L:34} Assume that $G$ has a somewhere dense
(resp. dense) orbit in $\mathbb{R}^{n}$. Then for every $v\in U$, $G(v)$ is somewhere dense
(resp. dense) in $\mathbb{R}^{n}$. In particular, $G(u_{\eta})$ is somewhere dense $($resp. dense$)$ in $\mathbb{R}^{n}$.
\end{lem}

\begin{proof} Let $u\in \mathbb{R}^{n}$ such that $\overset{\circ}{\overline{G(u)}}\neq\emptyset$. 
By Lemma \ref{L:00}, we have that $u\in U$. Let $v\in U$, then by 
Proposition \ref{pp:01},
 $v=Bu$ for some $B\in \textrm{Vect}(G)$. Moreover, by Lemma \ref{L:74}, (i), 
 $B\in \mathcal{K}^{*}_{\eta}(\mathbb{R})$. It follows that
$G(v)=B(G(u))$ and, since $B$ is invertible, 
$\overset{\circ}{\overline{G(v)}}\neq\emptyset$.

Now if $\overline{G(u)} = \mathbb{R}^{n}$, then $\overline{G(v)}= B\big(\overline{G(u)}\big)
= B\big(\mathbb{R}^{n}\big) = \mathbb{R}^{n}$.
\end{proof}

\begin{lem}\label{L:7+0} If \;$\overset{\circ}{\overline{G(u_{\eta})}}\neq\emptyset$ 
$($resp. $\overset{\circ}{\overline{\mathrm{g}_{\eta}(u_{\eta})}}\neq\emptyset$$)$, 
then \;$\mathcal{C}_{u_{\eta}} = \Phi_{u_{\eta}}\big(\mathcal{C}(G)\cap
\mathcal{K}_{\eta}^{+}(\mathbb{R})\big)$ and
\;$G^{+}(u_{\eta}) = G(u_{\eta})\cap \mathcal{C}_{u_{\eta}}$.
\end{lem}

\begin{proof} Let $\Phi_{u_{\eta}}(M)= Mu_{\eta}= v$ where $M\in \mathcal{C}(G)\cap\mathcal{K}^{+}_{\eta}(\mathbb{R})$.
Then $M = \textrm{diag}(M_{1},\dots, M_{r}, \widetilde{M_{1}},\dots,\widetilde{M_{s}})$, where
$M_{k}\in\mathbb{T}_{n_{k}}^{+}(\mathbb{R})$ and $\widetilde{M_{l}}\in\mathbb{B}_{m_{l}}^{*}(\mathbb{R})$ for every
$k=1,\dots, r$ and $l=1,\dots,s$. Write
$v=[v_{1},\dots,v_{r};\widetilde{v_{1}},\dots,\widetilde{v_{s}}]^{T}\in\mathbb{R}^{n}$, where
 $v_{k}\in\mathbb{R}^{n_{k}}$, $\widetilde{v_{l}}\in\mathbb{R}^{2m_{l}}$, $k=1,\dots, r$ and $l=1,\dots,s$. It follows that
 $M_{k}e_{\eta,k}=v_{k}$ and  $\widetilde{M_{l}}f_{\eta,l}=\widetilde{v_{l}}$.
Therefore $v_{k}\in \mathbb{R}^{+}\times \mathbb{R}^{n_{k}-1}$ and $\widetilde{v_{l}}\in (\mathbb{R}^{2}\backslash\{(0,0)\})\times
\mathbb{R}^{2m_{l}-2}$. We conclude that $v\in C_{u_{\eta}}$.

Conversely, let
$v=[v_{1},\dots,v_{r};\widetilde{v_{1}},\dots,\widetilde{v_{s}}]^{T}\in C_{u_{\eta}}$. Then $v_{k}\in \mathbb{R}^{+}\times \mathbb{R}^{n_{k}-1}$ and
$\widetilde{v_{l}}\in (\mathbb{R}^{2}\backslash\{(0,0)\})\times
\mathbb{R}^{2m_{l}-2}$. Since $\Phi_{u_{\eta}}$ is bijective (Proposition \ref{pp:01}),
there exists $B\in \mathcal{C}(G)$ such that $\Phi_{u_{\eta}}(B)= Bu_{\eta}= v$. By Lemma ~\ref{L:74},(ii), it follows that 
$B\in \mathcal{K}^{+}_{\eta}(\mathbb{R})$ and so $v\in \Phi_{u_{\eta}}\big(\mathcal{C}(G)\cap
\mathcal{K}_{\eta}^{+}(\mathbb{R})\big)$.

On the other hand, we have \begin{align*}
G(u_{\eta})\cap \mathcal{C}_{u_{\eta}}& = \Phi_{u_{\eta}}(G)\cap\Phi_{u_{\eta}}(\mathcal{C}(G)\cap\mathcal{K}^{+}_{\eta}(\mathbb{R}))\\
\ & \ = \Phi_{u_{\eta}}(G\cap\mathcal{K}^{+}_{\eta}(\mathbb{R})) = \Phi_{u_{\eta}}(G^{+}) = G^{+}(u_{\eta}).
\end{align*}
\end{proof}

\begin{lem}\label{L:007} Let $G$ be an abelian subgroup of
$\mathcal{K}^{*}_{\eta}(\mathbb{R})$. If $\overset{\circ}{\overline{G(u_{\eta})}}\neq\emptyset$, then for any
$v\in C_{u_{\eta}}$, we have $\overline{G(v)}\cap
C_{u_{\eta}} = C_{u_{\eta}}$. In particular, $G(u_{\eta})$ is dense in  $C_{u_{\eta}}$.
\end{lem}

\begin{proof}  Since $C_{u_{\eta}}$ is connected, it suffices to prove that \;$\overline{G(v)}\cap C_{u_{\eta}}
=\overset{\circ}{\overline{G(v)}}\cap C_{u_{\eta}}$: Suppose there exists
$w\in (\overline{G(v)} \ \backslash
\ \overset{\circ}{\overline{G(v)}})\cap C_{u_{\eta}}$. Since  $G$ is a group, $\overline{G(v)} \ \backslash
\ \overset{\circ}{\overline{G(v)}}$ is $G$-invariant and therefore
$\overline{G(w)}\subset \overline{G(v)} \ \backslash
\ \overset{\circ}{\overline{G(v)}}$. By Lemma ~\ref{L:34}, $\overset{\circ}{\overline{G(w)}}\neq\emptyset$.
Hence $ \overline{G(v)} \ \backslash
\ \overset{\circ}{\overline{G(v)}}$ has a non-empty interior, which is absurd. This completes the
proof.
\end{proof}

\begin{prop} \label{p:4} Assume that
$\overset{\circ}{\overline{G(u_{\eta})}}\neq\emptyset$  or
$\overset{\circ}{\overline{\mathrm{g}_{\eta}(u_{\eta})}}\neq\emptyset$, where $u_{\eta}$ is defined in paragraph 1. 
Then $f:=\Phi_{u_{\eta}}\circ
\mathrm{exp}_{|\mathcal{K}_{\eta}(\mathbb{R})}\circ \Phi_{u_{\eta}}^{-1}$
from $\mathbb{R}^{n}$ to $\mathbb{R}^{n}$ is well defined and
satisfies the following:
\begin{itemize}
    \item [(i)] $f$ is continuous and open;
  \item[(ii)]  $f(Bu_{\eta})=e^{B}u_{\eta}$ for every $B\in \mathcal{C}(G)$;
  \item [(iii)]$f^{-1}(G^{+}(u_{\eta}))= \mathrm{g}_{\eta}(u_{\eta})$ and
  $f(\mathrm{g}_{\eta}(u_{\eta}))=G^{+}(u_{\eta})$;
 \item [(iv)]  $f(\mathbb{R}^{n}) =  C_{u_{\eta}}$.
\end{itemize}
\end{prop}

\begin{proof} (i) By Proposition ~\ref{pp:01}, $\Phi_{u_{\eta}}$
 is a linear isomorphism. So $f:=\Phi_{u_{\eta}}\circ
exp_{|\mathcal{K}_{\eta}(\mathbb{R})}\circ \Phi_{u_{\eta}}^{-1}$
is well defined and continuous. Moreover, $f$ is a local
diffeomorphism by Corollary~\ref{LLL:1}, and therefore $f$ is an open
map.\ \\ (ii) For every $B\in \mathcal{C}(G)$, we have that
$\Phi_{u_{\eta}}^{-1}(Bu_{\eta})=B$. Therefore,
\begin{align*} f(Bu_{\eta})& =\Phi_{u_{\eta}}\circ
\textrm{exp}_{|\mathcal{K}_{\eta}(\mathbb{R})}\circ
\Phi_{u_{\eta}}^{-1}(Bu_{\eta})\\ \ & =\Phi_{u_{\eta}}(e^{B}) =e^{B}u_{\eta}.
\end{align*}

(iii) We have that  \begin{align*} f^{-1}(G^{+}(u_{\eta}))&
=\Phi_{u_{\eta}}\circ
\textrm{exp}^{-1}_{|\mathcal{K}_{\eta}(\mathbb{R})}\circ
\Phi_{u_{\eta}}^{-1}(G^{+}(u_{\eta}))\\ \ &
=\Phi_{u_{\eta}}(\textrm{exp}^{-1}_{|\mathcal{K}_{\eta}(\mathbb{R})}(G^{+}))
=\Phi_{u_{\eta}}(\mathrm{g}_{\eta}) =\mathrm{g}_{\eta}(u_{\eta}).
\end{align*}

 We also have that \begin{align*} f(\textrm{g}_{u_{\eta}})&
= \Phi_{u_{\eta}}\circ \textrm{exp}_{|\mathcal{K}_{\eta}(\mathbb{R})}\circ
\Phi_{u_{\eta}}^{-1}(\mathrm{g}_{u_{\eta}})\\
& = \Phi_{u_{\eta}}(\textrm{exp}_{|\mathcal{K}_{\eta}(\mathbb{R})}(\textrm{g}_{\eta}))
= \Phi_{u_{\eta}}(G^{+}) = G^{+}(u_{\eta}).
\end{align*}
\
\\
(iv) As $\Phi^{-1}_{u_{\eta}} : \mathbb{R}^{n}\longrightarrow
\mathcal{C}(G)$  is an isomorphism, we get that \begin{align*}
f(\mathbb{R}^{n})& = \Phi_{u_{\eta}}\circ
\textrm{exp}_{|\mathcal{K}_{\eta}(\mathbb{R})}\circ
\Phi_{u_{\eta}}^{-1}(\mathbb{R}^{n})\\
& = \Phi_{u_{\eta}}\left(\textrm{exp}_{|\mathcal{K}_{\eta}(\mathbb{R})}(\mathcal{C}(G))\right)\\
\ & =\Phi_{u_{\eta}}(\mathcal{C}(G)\cap
\mathcal{K}^{+}_{\eta}(\mathbb{R})) \qquad  (\textrm{by Lemma}~\ref{L:1},
(iii))\\ \ & = \mathcal{C}_{u_{\eta}} \qquad  \qquad \qquad  \qquad \; \; \;\;\;(\textrm{by Lemma}~\ref{L:7+0}).
\end{align*}
\end{proof}

\begin{prop}\label{p:5n} Let $G$ be an abelian sub-semigroup of
$\mathcal{K}^{*}_{\eta}(\mathbb{R})$. The following assertions are equivalent:
\smallskip
\begin{itemize}
  \item [(i)] $G(u_{\eta})$ is somewhere dense in $\mathbb{R}^{n}$,
  \item [(ii)] $\mathrm{g}_{\eta}(u_{\eta})$ is somewhere dense in $\mathbb{R}^{n}$.
\end{itemize}
\end{prop}

\begin{proof} $(ii) \Longrightarrow (i)$. Suppose that  
$\mathrm{g}_{\eta}(u_{\eta})$ is somewhere dense in $\mathbb{R}^{n}$ i.e. $\overset{\circ}{\overline{\textrm{g}_{\eta}(u_{\eta})}}\neq \emptyset$.  
By Proposition ~\ref{p:4}, there exists a continuous open map  $f :\mathbb{R}^{n}\longrightarrow \mathbb{R}^{n}$ such that
   $f(\textrm{g}_{\eta}(u_{\eta})) = G^{+}(u_{\eta})$. Hence, $f(\overset{\circ}{ \overline{\mathrm{g}_{\eta}(u_{\eta})})}\subset 
  \overset{\circ}{\overline{G^{+}(u_{\eta})}}$, and therefore
 $\overset{\circ}{\overline{G(u_{\eta})}}\neq\emptyset$.
 
 $(i)\Longrightarrow (ii)$. Suppose that  $\overset{\circ}{\overline{G(u_{\eta})}}\neq\emptyset$. Then 
 $f(\mathbb{R}^{n}) = C_{u_{\eta}}$ and $f^{-1}(G^{+}(u_{\eta})) =
\textrm{g}_{\eta}(u_{\eta})$, where $f$ is defined in Proposition ~\ref{p:4}. Moreover 
$\overset{\circ}{\overline{G(u_{\eta})}}\cap C_{u_{\eta}}\neq\emptyset$. Indeed, there is $A\in G$ such that $A(u_{\eta})\in \overset{\circ}{\overline{G(u_{\eta})}}$, 
so $A^{2}(u_{\eta})\in \overset{\circ}{\overline{G(u_{\eta})}}\cap C_{u_{\eta}}$. 

It follows that $\emptyset\neq f^{-1}(\overset{\circ}{\overline{G(u_{\eta})}}\cap C_{u_{\eta}})$. As $f$ is continuous and open and 
$G^{+}(u_{\eta}) = G(u_{\eta})\cap \mathcal{C}_{u_{\eta}}$ 
(Lemma \ref{L:7+0}), so 
$$f^{-1}(\overset{\circ}{\overline{G(u_{\eta})}}\cap C_{u_{\eta}})\subset f^{-1}(\overset{\circ}{\overline{G^{+}(u_{\eta})}}) \subset 
\overset{\circ}{\overline{f^{-1}(G^{+}(u_{\eta}))}} =
\overset{\circ}{\overline{\textrm{g}_{\eta}(u_{\eta})}}.$$ Hence \;
$\overset{\circ}{\overline{\textrm{g}_{\eta}(u_{\eta})}}\neq \emptyset$.
\end{proof}

\begin{prop}\label{p:5nn} Let $G$ be an abelian sub-semigroup of
$\mathcal{K}^{*}_{\eta}(\mathbb{R})$. The following assertions are equivalent:
\smallskip
\begin{itemize}
  \item [(i)] $G(u_{\eta})$ is dense in  $C_{u_{\eta}}$,
  \item [(ii)] $\mathrm{g}_{\eta}(u_{\eta})$ is dense in $\mathbb{R}^{n}$.
\end{itemize}
\end{prop}

\begin{proof} $(ii) \Longrightarrow (i)$. Suppose that $\mathrm{g}_{\eta}(u_{\eta})$ is dense in $\mathbb{R}^{n}$. 
By Proposition ~\ref{p:4}, there exists a continuous open map  $f :\mathbb{R}^{n}\longrightarrow \mathbb{R}^{n}$ such that
   $f(\textrm{g}_{\eta}(u_{\eta})) = G^{+}(u_{\eta})$. Hence, $C_{u_{\eta}} = f(\mathbb{R}^{n}) = f(\overline{\mathrm{g}_{\eta}(u_{\eta})}) \subset 
  \overline{G^{+}(u_{\eta})}$. As $G^{+}(u_{\eta}) = G(u_{\eta})\cap \mathcal{C}_{u_{\eta}}$ 
(Lemma \ref{L:7+0}), therefore $G(u_{\eta})$ is dense in  $C_{u_{\eta}}$.
  
  $(i)\Longrightarrow (ii)$. Suppose that $G(u_{\eta})$ is dense in $C_{u_{\eta}}$ i.e. $G^{+}(u_{\eta}) = G(u_{\eta})\cap \mathcal{C}_{u_{\eta}}$ is dense in  $C_{u_{\eta}}$. Then 
 $f(\mathbb{R}^{n}) = C_{u_{\eta}}$ and $f^{-1}(G^{+}(u_{\eta})) =
\textrm{g}_{\eta}(u_{\eta})$, where $f$ is defined in Proposition ~\ref{p:4}. As $f$ is continuous, so 
$$\mathbb{R}^{n} = f^{-1}(C_{u_{\eta}})\subset f^{-1}(\overline{G^{+}(u_{\eta})}) \subset 
\overline{f^{-1}(G^{+}(u_{\eta}))} =
\overline{\textrm{g}_{\eta}(u_{\eta})}.$$ Hence $\mathbb{R}^{n}  =
\overline{\textrm{g}_{\eta}(u_{\eta})}$.
 
\end{proof}

\begin{cor}\label{c:5} Let $G$ be an abelian subgroup of
$\mathcal{K}^{*}_{\eta}(\mathbb{R})$. The following assertions are equivalent:
\smallskip
\begin{itemize}
  \item [(i)] $G(u_{\eta})$ is somewhere dense in $\mathbb{R}^{n}$,
  \item [(ii)] $\mathrm{g}_{\eta}(u_{\eta})$ is dense in $\mathbb{R}^{n}$.
\end{itemize}
\end{cor}

\begin{proof} $(ii) \Longrightarrow (i)$ is already shown in Proposition \ref{p:5nn}.
 $(i)\Longrightarrow (ii)$. Suppose that  $\overset{\circ}{\overline{G(u_{\eta})}}\neq\emptyset$. Then
$\overline{G(u_{\eta})}\cap C_{u_{\eta}} = C_{u_{\eta}}$ (Lemma ~\ref{L:007}). Therefore $C_{u_{\eta}} = \overline{G(u_{\eta})\cap
C_{u_{\eta}}}\cap C_{u_{\eta}}$ i.e.  $G(u_{\eta})$ is dense in  $C_{u_{\eta}}$.
Hence by Proposition \ref{p:5nn}, \;
$\overline{\textrm{g}_{\eta}(u_{\eta})} = \mathbb{R}^{n}$.
\end{proof}

%\begin{prop}\label{p:1nn} Let $G$ be an abelian sub-semigroup of $M_{n}(\mathbb{R})$ and $P\in \textrm{GL}(n, \mathbb{R})$ such that $P^{-1}GP\subset
%\mathcal{K}_{\eta}(\mathbb{R})$ where $\eta$ has length $(r+2s)$. 
%\begin{enumerate}
% \item The following properties are equivalent:
% \begin{itemize}
%  \item [(i)] $G$ has a dense orbit in $P(C_{u_{\eta}})$,
%  \item [(ii)] The orbit $G(Pu_{\eta})$ is dense in $P(C_{u_{\eta}})$,
%  \item [(iii)] $\mathrm{g}_{\eta}(Pu_{\eta})$ is an additive sub-semigroup dense in $\mathbb{R}^{n}$.
%\end{itemize}

%\item Assume that $G$ is generated by $p$ matrices $A_{1},\dots,A_{p}$  ($p\geq 1$)
%and let $B_{1},\dots,B_{p}\in \mathrm{g}_{\eta}$ such that $A_{1}^{2} =
%e^{B_{1}},\dots, A_{p}^{2} = e^{B_{p}}$. Then $G$ has dense orbit in  $P(C_{u_{\eta}})$ if and only if
%$$\underset{k=1}{\overset{p}{\sum}}\mathbb{N}B_{k}Pu_{\eta} +
%\underset{l=1}{\overset{s}{\sum}}2\pi\mathbb{Z}Pf_{\eta}^{(l)}$$ is dense in $\mathbb{R}^{n}$.
%\end{enumerate}
%\end{prop}
%\begin{proof}
%\end{proof}

\section{\bf Proof of the first part of Theorems \ref{T:1n}, ~\ref{T:2} and Corollary \ref{c:1}}

\subsection{Proof of the first part of Theorem ~\ref{T:1n} and Corollary \ref{c:1}}

\begin{prop}\label{p:6} Let $G$ be an abelian sub-semigroup of $M_{n}(\mathbb{R})$ and let $u\in \mathbb{R}^{n}$.
Then $G^{*}(u)$ is  somewhere dense $($resp. dense$)$ in $\mathbb{R}^{n}$ if and only if $G(u)$ is.
\end{prop}

\begin{proof} The first implication is trivial. Conversely, suppose that
$\overset{\circ}{\overline{G(u)}}\neq \emptyset$ (resp.
$\overline{G(u)}=\mathbb{R}^{n}$). We can assume, using
Proposition ~\ref{p:2}, that
 $G\subset\mathcal{K}_{\eta}(\mathbb{R})$. We let $G^{\prime}:=G\backslash G^{*}$.
\
 \\
 - If $G^{\prime} = \emptyset$, then $G=G^{*}$ and so $\overset{\circ}{\overline{G^{*}(u)}}\neq \emptyset$ (resp. $\overline{G^{*}(u)} =  \mathbb{R}^{n}$).
\
 \\
 - If $G^{\prime}\neq \emptyset$, then

$$G(u)\subset\left(\underset{A\in G^{\prime}
}{\bigcup}\textrm{Im}(A)\right)\cup G^{*}(u).$$
\
\\
Since every $A\in G^{\prime}$ is non invertible,
$\textrm{Im}(A)\subset \underset{k=1
}{\overset{r}{\bigcup}}H_{k}\cup \underset{l=1
}{\overset{s}{\bigcup}}F_{l}$, where $$H_{k}:= \left\{u = [u_{1},\dots,u_{r};\
\widetilde{u}_{1},\dots,\widetilde{u}_{s}]^{T}\in \mathbb{R}^{n}: \  u_{j}\in \mathbb{R}^{n_{j}},
   u_{k}\in\{0\}\times\mathbb{R}^{n_{k}-1}, \ 1\leq j\neq k\leq r \right\}$$
and $$F_{l}:=\left\{u=[u_{1},\dots,u_{r};\
\widetilde{u}_{1},\dots,\widetilde{u}_{s}]^{T}\in \mathbb{R}^{n}: \
 \widetilde{u}_{l}\in\{(0,0)\}\times\mathbb{R}^{2m_{l}-2}, \ 1\leq l\leq s\right\}.$$

\
\\
It follows that
  $$G(u)\subset\left(\underset{k=1 }{\overset{r}{\bigcup}}H_{k}\cup \underset{l=1 }{\overset{s}{\bigcup}}F_{l}\right)\cup G^{*}(u)$$ and so
$$\overline{G(u)}\subset\left(\underset{k=1
}{\overset{r}{\bigcup}}H_{k}\cup \underset{l=1
}{\overset{s}{\bigcup}}F_{l}\right)\cup \overline{G^{*}(u)}.$$
\
\\
 Since  $H_{k}$ (resp. $F_{l}$) has dimension $n-1$ (resp. $n-2$), $\overset{\circ}{H_{k}}= 
 \overset{\circ}{F_{l}} = \emptyset$, for every $1\leq k\leq r$, $1\leq l \leq s$, and therefore
$\overset{\circ}{\overline{G^{*}(u)}}\neq \emptyset$ (resp.
$\overline{G^{*}(u)} = \mathbb{R}^{n}$).$\mathrm{g}_{\eta}(u_{\eta})$ is dense in $\mathbb{R}^{n}$.
\end{proof}
\smallskip

\begin{lem}\label{L:6} Let $G$ be an abelian sub-semigroup of $\mathcal{K}_{\eta}(\mathbb{R})$ and
 $\mathrm{g}_{\eta}^{*} = \mathrm{exp}^{-1}(G^{*})\cap \mathcal{K}_{\eta}(\mathbb{R})$. Then
$\mathrm{g}_{\eta} = \mathrm{g}_{\eta}^{*}$.
\end{lem}

\begin{proof} Since $G^{*}\subset G$, we see that $\textrm{g}_{\eta}^{*}\subset \textrm{g}_{\eta}$.
Conversely, if $B\in\textrm{g}_{\eta}$, then $e^{B}\in G\cap \textrm{GL}(n,
\mathbb{R})=G^{*}$. So $B\in \textrm{exp}^{-1}(G^{*})\cap
\mathcal{K}_{\eta}(\mathbb{R})=\textrm{g}_{\eta}^{*}$ and hence
 $\textrm{g}_{\eta}\subset \textrm{g}_{\eta}^{*}$. It follows that $\textrm{g}_{\eta} = \textrm{g}_{\eta}^{*}$.
\end{proof}
\smallskip

\begin{cor}\label{c:38n} If all matrices of $G\backslash I_{n}$ are non
invertible, then all orbits of $G$ are
nowhere dense.
\end{cor}

\begin{proof} This follows from  Proposition ~\ref{p:6} since $G^{*} = \{I_{n}\}$. 
\end{proof}
\medskip

\textit{Proof of the first part of Theorem \ref{T:1n} and Corollary \ref{c:1}}.  One can assume by
Proposition ~\ref{p:2} that $G$ is an abelian sub-semigroup of
$\mathcal{K}_{\eta}(\mathbb{R})$. From Proposition \ref{p:6} and Lemma \ref{L:6}, we may assume that 
$G\subset \mathcal{K}_{\eta}^{*}(\mathbb{R})$.

$(ii)\Longrightarrow (i)$. This is clear.  

$(i)\Longrightarrow (ii)$.  This follows directly from Lemma ~\ref{L:34}. 

The equivalence $(ii)\Longleftrightarrow (iii)$ in Theorem \ref{T:1n} follows from Proposition ~\ref{p:5n} and that in Corollary \ref{c:1}, follows from 
Corollary ~\ref{c:5} if $G$ is 
a group and from Corollary \ref{c:7n} if $G$ is a finitely generated abelian sub-semigroup. \qed
\medskip

%\textit{Proof of the first part of Theorem \ref{T:1}}. 
%
%$(ii)\Longrightarrow (i)$. This is clear.  
%
%$(i)\Longrightarrow (ii)$.  This follows directly from Lemma ~\ref{L:34}. 
%
%$(ii)\Longleftrightarrow (iii)$. This follows from Corollary ~\ref{c:5}. \qed
%\medskip

\subsection{ Proof of the first part of Theorem \ref{T:2}}
\medskip

Denote by $\Gamma$ \;the subgroup of
$\mathcal{K}^{*}_{\eta}(\mathbb{R})$ generated by $(S_{k})_{1\leq k\leq r}$, where 
\
\\

$S_{k}: =
\textrm{diag}\left(\varepsilon_{1,k}I_{n_{1}},\dots,
\varepsilon_{r,k}I_{n_{r}}; \ I_{2m_{1}} ,\dots,
I_{2m_{s}}\right)\in \mathcal{K}^{*}_{\eta}(\mathbb{R}),$
and 
\
\\
$$\varepsilon_{i,k}:= \begin{cases}

        -1 & \ \textrm{if} \ \ {i=k}, \\
  1 & \ \textrm{if} \ \ {i\neq k}. \end{cases}, \ \ 1\leq i, \ k\leq r$$

\begin{lem}\label{L:7+} Let  $G$  be an abelian sub-semigroup of
$\mathcal{K}^{*}_{\eta}(\mathbb{R})$. Then:
\medskip

\begin{itemize}

\item[(i)] $U = \underset{S\in\Gamma}{\bigcup}S(C_{u_{\eta}})$.

\item[(ii)] $S_{k}M = MS_{k}$, for every $M\in \mathcal{K}^{*}_{\eta}(\mathbb{R})$,  $k=1,\dots,r$.

\item[(iii)] if \;$\mathrm{ind}(G) = r$, then $G(u_{\eta})\cap S(C_{u_{\eta}})\neq\emptyset$ for every  $S\in\Gamma$.
 \end{itemize}
\end{lem}

\begin{proof} Assertions (i) and (ii) are easier to prove. Assertion (iii). By hypothesis and for every 
$1\leq k\leq r$, there exists $B^{(k)}\in G$ such that $B^{(k)}_{k}$ has only one eigenvalue
$\mu_{k}<0$ and all its other real eigenvalues $\mu_{i}>0$, $i\neq k$. As 
$B^{(k)}\mathcal{C}(G)=S_{k}\mathcal{C}(G) = \mathcal{C}(G)$, then $B^{(k)}\Phi_{u_{\eta}}(\mathcal{C}(G)) = 
\Phi_{u_{\eta}}\big(B^{(k)}\mathcal{C}(G)\big)$, \\ ~ $S_{k}\Phi_{u_{\eta}}(\mathcal{C}(G)) = \Phi_{u_{\eta}}(S_{k}\mathcal{C}(G))$
and $B^{(k)}\mathcal{K}^{+}_{\eta}(\mathbb{R}) = S_{k}\mathcal{K}^{+}_{\eta}(\mathbb{R})$. So
\begin{align*}
B^{(k)}(C_{u_{\eta}}) & = B^{(k)}\Phi_{u_{\eta}}(\mathcal{C}(G)\cap\mathcal{K}^{+}_{\eta}(\mathbb{R}))\\
\ & = \Phi_{u_{\eta}}(B^{(k)}\big(\mathcal{C}(G)\cap\mathcal{K}^{+}_{\eta}(\mathbb{R}))\big)\\
\ & = \Phi_{u_{\eta}}(S_{k}\mathcal{C}(G)\cap S_{k}\mathcal{K}^{+}_{\eta}(\mathbb{R}))\\
\ & = S_{k}\Phi_{u_{\eta}}(\mathcal{C}(G)\cap \mathcal{K}^{+}_{\eta}(\mathbb{R}))\\
\ & = S_{k}C_{u_{\eta}}.
\end{align*} Thus $B^{(k)}(u_{\eta})\in S_{k}(C_{u_{\eta}})$ and hence $G(u_{\eta})\cap S_{k}(C_{u_{\eta}})\neq\emptyset$ for every $1\leq k \leq r$.
It follows that G($u_{\eta})\cap S(C_{u_{\eta}})\neq\emptyset$, for every $S\in\Gamma$. This completes the proof.
\end{proof}
\medskip

\begin{prop}\label{p:7+n}  Let  $G$ be an abelian sub-semigroup of
$\mathcal{K}^{*}_{\eta}(\mathbb{R})$. Then the following properties are equivalent:
\medskip

\begin{itemize}
\item [(i)] ~$\overline{G(u_{\eta})} = \mathbb{R}^{n}$

\item [(ii)] ~ $\mathrm{g}_{\eta}(u_{\eta})$ is dense in $\mathbb{R}^{n}$ and $\mathrm{ ind}(G) = r$.
\end{itemize}
\end{prop}
\medskip

\begin{proof} $(i)\Longrightarrow (ii)$. Assume that ~$\overline{G(u_{\eta})} = \mathbb{R}^{n}$. The fact that $\mathrm{g}_{\eta}(u_{\eta})$ is dense in $\mathbb{R}^{n}$ 
follows from Proposition \ref{p:5nn}. Now suppose that 
$\textrm{ind}(G) < r$.
Then there exists  $1\leq k_{0}\leq r$  such that for every
 $B = \textrm{diag}(B_{1},\dots,B_{r};\ \widetilde{B}_{1},\dots,\widetilde{B}_{s})\in G$ with
$B_{k}\in \mathbb{T}_{n_{k}}(\mathbb{R})$, $k = 1,\dots, r$ having
an eigenvalue $\mu_{k}$ and $\widetilde{B}_{l}\in
\mathbb{B}_{m_{l}}^{*}(\mathbb{R})$, $l = 1,\dots, s$, we have
$\mu_{k_{0}}>0$ or $\mu_{i}<0$ for some $i\neq k_{0}$. Therefore
$G(u_{\eta})\subset \mathbb{R}^{n}\backslash
\mathcal{C}^{\prime}_{u_{\eta},k_{0}}$, where
$\mathcal{C}^{\prime}_{u_{\eta},k_{0}}:=
S_{k_{0}}(\mathcal{C}_{u_{\eta}})$ and thus
$\overline{\mathbb{R}^{n}\backslash
\mathcal{C}^{\prime}_{u_{\eta},k_{0}}}= \mathbb{R}^{n}$; that is
$\mathcal{C}^{\prime}_{u_{\eta},k_{0}}=\emptyset$, a contradiction.\\

$(ii)\Longrightarrow (i)$. Suppose that  $\mathrm{g}_{\eta}(u_{\eta})$ is dense in $\mathbb{R}^{n}$ and $\mathrm{ ind}(G) = r$. 
By Lemma ~\ref{L:7+},(iii), we have
$G(u_{\eta})\cap S(C_{u_{\eta}})\neq\emptyset,$  for every $S\in\Gamma$.
So let  $B\in G$ such that $B(u_{\eta})\in S(C_{u_{\eta}})$. By Proposition \ref{p:5nn}, $\overline{G(u_{\eta})}\cap C_{u_{\eta}} = C_{u_{\eta}}$.
Moreover $B(C_{u_{\eta}}) = S(C_{u_{\eta}})$ (since  $B(C_{u_{\eta}})$ is a connected component of $U$ meeting $S(C_{u_{\eta}})$).
It follows that $$ B\left(\overline{G(u_{\eta})}\cap C_{u_{\eta}}\right) = \overline{BG(u_{\eta})}\cap
S(C_{u_{\eta}})= S(C_{u_{\eta}})$$ and hence $S(C_{u_{\eta}})\subset \overline{G(u_{\eta})}$,  for every $S\in\Gamma$. As $U =
\underset{S\in\Gamma}{\bigcup}S(C_{u_{\eta}})$ (Lemma ~\ref{L:7+},
(i)), then  $U\subset \overline{G(u_{\eta})}$
 and since $\overline{U} = \mathbb{R}^{n}$, it follows that $\overline{G(u_{\eta})} = \mathbb{R}^{n}$.
\end{proof}

\begin{prop}\label{p:7+}  Let $G$ be an abelian subgroup of
$\mathcal{K}^{*}_{\eta}(\mathbb{R})$. Then the following properties are equivalent:
\medskip

\begin{itemize}
\item [(i)] ~$\overline{G(u_{\eta})} = \mathbb{R}^{n}$

\item [(ii)] ~$\overset{\circ}{\overline{G(u_{\eta})}}\neq \emptyset$ $\mathrm{and}$ $\mathrm{ ind}(G) = r$.
\end{itemize}
\end{prop}

\begin{proof}  $(i)\Longrightarrow (ii)$ is already shown in Proposition \ref{p:7+n}. \\
$(ii)\Longrightarrow (i)$. Suppose that  ~$\overset{\circ}{\overline{G(u_{\eta})}}\neq \emptyset$ and $\mathrm{ ind}(G) = r$. 
By Lemma ~\ref{L:7+},(iii), we have
$G(u_{\eta})\cap S(C_{u_{\eta}})\neq\emptyset,$  for every $S\in\Gamma$.
So let  $v\in G(u_{\eta})\cap S(C_{u_{\eta}})$ and
 $w = S^{-1}(v)\in C_{u_{\eta}}$. By Lemma ~\ref{L:007}, $\overline{G(w)}\cap C_{u_{\eta}} = C_{u_{\eta}}$ 
 and by Lemma ~\ref{L:7+}, (ii),   $G(v)= G(Sw) = S(G(w))$. It follows that $$\overline{G(v)}\cap
S(C_{u_{\eta}}) = S\left(\overline{G(w)}\cap C_{u_{\eta}}\right) =
S(C_{u_{\eta}}),$$ and hence $S(C_{u_{\eta}})\subset \overline{G(u_{\eta})}$. As $U =
\underset{S\in\Gamma}{\bigcup}S(C_{u_{\eta}})$ (Lemma ~\ref{L:7+},
(i)), then  $U\subset \overline{G(u_{\eta})}$
 and since $\overline{U} = \mathbb{R}^{n}$, it follows that $\overline{G(u_{\eta})} = \mathbb{R}^{n}$.
\end{proof}

\begin{proof}[\it Proof of the first part of Theorem ~\ref{T:2}]
  We may assume that 
$G\subset \mathcal{K}_{\eta}^{*}(\mathbb{R})$.
$(i)\Longleftrightarrow (ii)$. This follows directly from Lemma ~\ref{L:34}.
 $(ii)\Longleftrightarrow (iii)$ follows from Proposition ~\ref{p:7+n}. 
 %(resp. Corollary \ref{c:5} and Proposition ~\ref{p:7+}). 
\end{proof}
\smallskip

\section{\bf Finitely generated abelian semigroups of $\mathcal{K}^{*}_{\eta}(\mathbb{R})$
with a somewhere dense (resp. dense) orbit}

 \begin{lem}\label{L:7}$($\cite{aAhM3}, Proposition 2.7$)$\ \begin{itemize}
 \item [(i)] Let $A,B\in \mathbb{T}_{n}(\mathbb{R})$. If  $e^{A} = e^{B}$, then  $A = B$.
\item [(ii)] Let  $A$,  $B\in \mathbb{B}_{m}(\mathbb{R})$. If  $e^{A} = e^{B}$, then $A = B + 2k\pi J_{m}$ for
some $k\in \mathbb{Z}$, where $J_{m} =\mathrm{diag}(J_{2},\dots,J_{2})\in \mathrm{GL}(2m, \ \mathbb{R})
 \ \ \mathrm{with} \ \ J_{2}=
  \left[\begin{array}{cc}
    0 & -1 \\
    1 & 0 \\
  \end{array}\right]$.
  \end{itemize}
\end{lem}

\begin{prop}\label{p:8} Let $G$ be an abelian sub-semigroup of $\mathcal{K}^{+}_{\eta}(\mathbb{R})$ and let
$B_{1},\dots,B_{p} \in  \mathcal{K}_{\eta}(\mathbb{R})$ $(p\geq 1)$ be such
that $e^{B_{1}},\dots, e^{B_{p}}$
 generate $G$. We have that $$\mathrm{g}_{u_{\eta}} = \begin{cases}
 \underset{k=1}{\overset{p}{\sum}}
 \mathbb{N}B_{k}u_{\eta} + \underset{l=1}{\overset{s}{\sum}} 2\pi \mathbb{Z}f_{\eta}^{(l)} & \ \mathrm{if } \ s\neq 0\\
  \underset{k=1}{\overset{p}{\sum}}
 \mathbb{N}B_{k}u_{\eta} & \ \mathrm{if } \ s=0
\end{cases}$$
 \end{prop}

\begin{proof} {\it Case $s\neq 0$}. 
\
\\
$\bullet$ First we determine $\textrm{g}$. Let $C\in \textrm{g}_{\eta}$. Then $C = \textrm{diag}(C_{1},\dots, C_{r};\
\widetilde{C}_{1},\dots, \widetilde{C}_{s} )\in
\mathcal{K}_{n}(\mathbb{R})$ and $e^{C} \in G$. So $e^{C} = e^{m_{1}B_{1}}\dots e^{m_{p}B_{p}}$ for some
$m_{1},\dots,m_{p}\in \mathbb{N}$. Since $B_{1},\dots,B_{p} \in
\textrm{g}_{\eta}$, they pairwise commute (Lemma ~\ref{L:1},(i)).
Therefore, $e^{C} = e^{m_{1}B_{1}+\dots+m_{p}B_{p}}$. Write
$B_{j} = \textrm{diag}(B_{j,1},\dots, B_{j,r};\
\widetilde{B}_{j,1},\dots, \widetilde{B}_{j,s})$. Then $e^{C_{k}}
= e^{m_{1}B_{1,k}+\dots+m_{p}B_{p,k}}$ and
$e^{\widetilde{C}_{l}} =
e^{m_{1}\widetilde{B}_{1,l}+\dots+m_{p}\widetilde{B}_{p,l}}$, $k = 1,\dots,r$, $l =
1,\dots,s$. As $C\in \textrm{g}_{\eta}$, we have $CB_{j} = B_{j}C$, and so $C_{k}B_{j,k} = B_{j,k}C_{k}$, 
and $\widetilde{C}_{l}\widetilde{B}_{j,l} = \widetilde{B}_{j,l}\widetilde{C}_{l},$
$k=1, \dots, r, \ j = 1,\dots,p; \ l=1, \dots, s$. From Lemma ~\ref{L:7}, it follows that
$C_{k} = m_{1}B_{1,k}+\dots +m_{p}B_{p,k}$ and $\widetilde{C}_{l}
= m_{1}\widetilde{B}_{1,l}+\dots +m_{p}\widetilde{B}_{p,l} + 2\pi
t_{l}J_{m_{l}}$ for some  $t_{l}\in \mathbb{Z}$, where $J_{m_{l}}
=\textrm{diag}(J_{2},\dots,J_{2})\in \textrm{GL}(2m_{l}, \
\mathbb{R})$.
 Therefore,
\medskip

\begin{align*}
C
&=\textrm{diag}\left(\underset{j=1}{\overset{p}{\sum}}m_{j}B_{j,1},\dots,
\underset{j=1}{\overset{p}{\sum}}m_{j}B_{j,r};
\underset{j=1}{\overset{p}{\sum}}m_{j} \widetilde{B}_{j,1}+2\pi
t_{1}J_{m_{1}},\ \dots,
\underset{j=1}{\overset{p}{\sum}}m_{j}\widetilde{B}_{j,s}+2\pi
t_{s}J_{m_{s}}\right)\\ \ &
=\underset{j=1}{\overset{p}{\sum}}m_{j}B_{j} +
\textrm{diag}\left(0,\dots,0;\ 2\pi t_{1}J_{m_{1}},\dots,2\pi
t_{s}J_{m_{s}}\right).
\end{align*}
\medskip
\
\\
Set $$L_{l} := \textrm{diag}(0,\dots,0;\
\widetilde{L}_{l,1},\dots, \widetilde{L}_{l,s})$$ with
$$\widetilde{L}_{l,i} =\left\{\begin{array}{c}
                   0 \in B_{m_{i}}(\mathbb{R}), \ \ \ \textrm{if} \ \ i \neq l, \\
                   J_{m_{l}}, \ \ \ \ \ \ \ \ \ \ \ \ \ \ \textrm{if} \ \   i = l. \
                 \end{array}\right.$$
\medskip
\
\\
Then we have $$\textrm{diag}(0,\dots, 0,\  2\pi
t_{1}J_{m_{1}},\dots, 2\pi t_{s}J_{m_{s}})
=\underset{l=1}{\overset{s}{\sum}}2\pi t_{l}L_{l},$$ and therefore $C=\underset{j=1}{\overset{p}{\sum}}m_{j}B_{j}+\underset{l=1}{\overset{s}{\sum}}2\pi
t_{l}L_{l}$. We conclude that
$$\textrm{g}_{\eta}=\underset{j=1}{\overset{p}{\sum}}\mathbb{N}B_{j}+\underset{l=1}{\overset{s}{\sum}}2\pi
\mathbb{Z}L_{l}.$$
\
\\
$\bullet$ \ Second, we determine $\textrm{g}_{\eta}(u_{\eta})$. Let $B\in
\textrm{g}_{\eta}$. We have
$$B=\underset{j=1}{\overset{p}{\sum}}m_{j}B_{j}+\underset{l=1}{\overset{s}{\sum}}2\pi
t_{l}L_{l}$$ for some $m_{1},\dots,m_{p}\in \mathbb{N}$, and
$t_{1},\dots, t_{s}\in \mathbb{Z}$. As
$\widetilde{L}_{l,i}f_{i,1}= f^{(l)}_{i}$, $i=1, \dots, s$, then
\begin{align*}
L_{l}u_{\eta} & = \textrm{diag}(0,\dots,\widetilde{L}_{l,1},\dots,
\widetilde{L}_{l,s})[e_{\eta,1},\dots,e_{\eta,r};\  f_{\eta,1}, \dots,
f_{\eta,s}]^{T}\\ \ & = [0,\dots, 0;\ f^{(l)}_{1},\dots,
f^{(l)}_{s}]^{T}\\ \ & = f_{\eta}^{(l)}.
\end{align*}

Hence, $Bu_{\eta} = \underset{j=1}{\overset{p}{\sum}}m_{j}B_{j}u_{\eta}+\underset{l=1}{\overset{s}{\sum}}2\pi
t_{l}f_{\eta}^{(l)}$, and therefore $$\textrm{g}_{\eta}(u_{\eta}) =
\underset{j=1}{\overset{p}{\sum}}\mathbb{N}B_{j}u_{\eta}+\underset{l=1}{\overset{s}{\sum}}2\pi\mathbb{Z}f_{\eta}^{(l)}.$$

{\it Case $s= 0$}. The same proof as before works; we get
$\textrm{g}_{\eta}=\underset{j=1}{\overset{p}{\sum}}\mathbb{N}B_{j}$, and so $\textrm{g}_{\eta}
(u_{\eta}) =
\underset{j=1}{\overset{p}{\sum}}\mathbb{N}B_{j}u_{\eta}$. This proves the proposition.
\end{proof}
\bigskip

Let $G$ be an abelian sub-semigroup of $\mathcal{K}_{\eta}^{*}(\mathbb{R})$. Denote by:

\
\\
\textbullet \;$G^{2} = \{A^{2}: A\in G\}$
\
\\
\textbullet \;$\textrm{g}_{\eta}^{2} = \textrm{exp}^{-1}(G^{2})\cap \left[
P(\mathcal{K}_{\eta}(\mathbb{R}))P^{-1}\right]$
\
\\
\textbullet \;$G^{*2} = \{A^{2} :\ A\in G^{*}\}$
\medskip

\begin{lem}\label{L:8} Let  $G$ be an abelian sub-semigroup of $\mathcal{K}_{\eta}^{*}(\mathbb{R})$.
Then: $$\overset{\circ}{\overline{G(u_{\eta})}}\neq \emptyset\ \ \
\mathrm{ if \ and \ only \ if } \ \
\overset{\circ}{\overline{G^{2}(u_{\eta})}}\neq
 \emptyset.$$
\end{lem}

\begin{proof} The `` if part'' is trivial. For the `` only if part'': suppose that
$\overset{\circ}{\overline{G(u_{\eta})}}\neq \emptyset$. Then by
Theorem ~\ref{T:1n}, $\overset{\circ}{\overline{\textrm{g}_{\eta}(u_{\eta})}}\neq \emptyset$. As $\textrm{g}_{\eta}\subset \frac{1}{2}\textrm{g}_{\eta}^{2}$
(since if $B\in \textrm{g}_{\eta}$, we have $e^{2B}= (e^{B})^{2}\in
G^{2}$), then
 $\overset{\circ}{\overline{\frac{1}{2}\textrm{g}_{\eta}^{2}(u_{\eta})}}\neq \emptyset$, and so  
 $\overset{\circ}{\overline{\textrm{g}_{\eta}^{2}(u_{\eta})}}\neq \emptyset$. Applying
Theorem ~\ref{T:1n} to the abelian sub-semigroup $G^{2}$, it
follows that
$\overset{\circ}{\overline{G^{2}(u_{\eta})}}\neq\emptyset$.
\end{proof}

\begin{cor}\label{C:45} Let $G$ be an abelian sub-semigroup of $\mathcal{K}^{*}_{\eta}(\mathbb{R})$. Then $G$ has
 a somewhere dense orbit if and only if so does $G^{2}$.
\end{cor}

\begin{proof} This is a consequence from Lemmas \ref{L:34} and \ref{L:8}.
\end{proof}

 \begin{lem}\label{l:31n} $\mathrm{g}_{\eta}(u_{\eta})$ is an additive sub-semigroup dense (resp. somewhere dense) in 
 $\mathbb{R}^{n}$ if and only if ~$(\mathrm{g}_{\eta})^{2}(u_{\eta})$ is.
 \end{lem}
 
 \begin{proof} This follows from the fact that $\mathrm{g}_{\eta}\subset \frac{1}{2}(\mathrm{g}_{\eta})^{2}$ and 
 $(\mathrm{g}_{\eta})^{2}\subset \mathrm{g}_{\eta}$.
 \end{proof}
\bigskip

\begin{proof}[Proof of the second part of Theorem \ref{T:1n}] \
This follows from the first part of Theorem
~\ref{T:1n}, Lemma \ref{L:8}, Proposition ~\ref{p:8} and the fact that 
an additive finitely generated semigroup that is somewhere dense in
$\mathbb{R}^{n}$ is dense in $\mathbb{R}^{n}$ (cf. \cite{hAaM}, Theorem 2.1).
\end{proof}

\begin{proof}[Proof of the second part of Corollary \ref{c:1}] \
This follows from the first part of Corollary \ref{c:1}, Lemma \ref{L:8} and that 
$\textrm{g}_{\eta}^{2}(u_{\eta}) = \underset{k=1}{\overset{p}{\sum}}\mathbb{Z}B_{k}Pu_{\eta} +
\underset{l=1}{\overset{s}{\sum}}2\pi\mathbb{Z}Pf_{\eta}^{(l)}$ by Proposition ~\ref{p:8}.
\end{proof}

\begin{proof}[Proof of the second part of Theorem \ref{T:2}] \
This follows from the second part of Theorem \ref{T:1n}, the first part of Theorem
~\ref{T:2}, Lemmas \ref{l:31n} and \ref{L:6} and Proposition ~\ref{p:8}.
\end{proof}

\begin{cor}\label{c:7n} Let $G$ be a finitely generated
abelian sub-semigroup of
$\mathcal{K}^{*}_{\eta}(\mathbb{R})$. The following assertions are equivalent:
\smallskip
\begin{itemize}
  \item [(i)] $G(u_{\eta})$ is somewhere dense in $\mathbb{R}^{n}$,
  \item [(ii)] $G(u_{\eta})$ is dense in  $C_{u_{\eta}}$,
  \item [(iii)] $\mathrm{g}_{\eta}(u_{\eta})$ is dense in $\mathbb{R}^{n}$.
\end{itemize}
\end{cor}

\begin{proof} $(ii) \Longleftrightarrow (iii)$ is already shown in Proposition \ref{p:5nn}.\\
 $(i)\Longrightarrow (iii)$. Suppose that $G(u_{\eta})$ is somewhere dense in $\mathbb{R}^{n}$. Then by Proposition \ref{c:5}, $\mathrm{g}_{\eta}(u_{\eta})$ is somewhere dense in 
$\mathbb{R}^{n}$. By Proposition \ref{p:8}, $\textrm{g}_{\eta}^{2}(u_{\eta})$ is finitely generated sub-semigroup of
$\mathbb{R}^{n}$ and thus so is $2\mathrm{g}_{\eta}(u_{\eta})$ (since $\mathrm{g}_{\eta}\subset \frac{1}{2}(\mathrm{g}_{\eta})^{2}$ 
and $(\mathrm{g}_{\eta})^{2}\subset \mathrm{g}_{\eta}$). By 
(\cite{hAaM}, Theorem 2.1), $2\mathrm{g}_{\eta}(u_{\eta})$ is dense in $\mathbb{R}^{n}$ and thus so is $\mathrm{g}_{\eta}(u_{\eta})$.
\end{proof}

\begin{cor}\label{c} If $G$ is a finitely generated
abelian sub-semigroup of
$\mathcal{K}^{*}_{\eta}(\mathbb{R})$, then Proposition \ref{p:7+} holds.
%Then the following
%properties are equivalent:
%\medskip
%
%\begin{itemize}
%\item [(i)] ~$\overline{G(u_{\eta})} = \mathbb{R}^{n}$
%
%\item [(ii)] ~$\overset{\circ}{\overline{G(u_{\eta})}}\neq \emptyset \ \mathrm{ and } \ \mathrm{ ind}(G) = r$.
%\end{itemize}
\end{cor}

\begin{proof}  $(i)\Longrightarrow (ii)$ is already shown in Proposition \ref{p:7+n}. 
$(ii)\Longrightarrow (i)$. This follows from Corollary \ref{c:7n} and Proposition \ref{p:7+n}.
%(or also from the second part of Theorems \ref{T:1n} and \ref{T:2}).
\end{proof}

\begin{prop}\label{p:9}$($\cite{sh}, Lemma 2.1$)$. Let $H = \mathbb{Z}u_{1}+\dots+\mathbb{Z}u_{m}$
with $u_{k}\in\mathbb{R}^{n}$,
$k = 1,\dots, m$. If $m\leq n$, then $H$ is nowhere dense in
$\mathbb{R}^{n}$.
\end{prop}
\
\\
{\it Proof of Corollary ~\ref{C:3}.} By Proposition \ref{p:9},
$\underset{k=1}{\overset{n-s}{\sum}}\mathbb{N}B_{k}Pu_{\eta} +
\underset{l=1}{\overset{s}{\sum}}2\pi\mathbb{Z}Pf_{\eta}^{(l)}$ is
nowhere dense in $\mathbb{R}^{n}$, and by applying
Theorem ~\ref{T:1n}, $G$ has nowhere dense orbit in
$\mathbb{R}^{n}$. \qed

\
\\
{\it Proof of Corollary ~\ref{C:4}}. Since $r+2s\leq n$, we see that $1+2s\leq n$ if $r\geq 1$ 
and $s\leq \frac{n}{2}$ if $r=0$. Hence $\left[\frac{n+1}{2}\right]\leq n-s$. Therefore the Corollary \ref{C:4} 
follows from Corollary \ref{C:3}. \qed
 \medskip

\section{\bf Proof of Theorem \ref{T:4} and Corollary \ref{C:6n}}

 We construct for every $n\in\mathbb{N}$, $n\geq 1$,
$r,s\in \mathbb{N}$, and for every partition
$\eta$ of $n$ of length $(r+2s)$, \; $(n-s+1)$ matrices
$A_{1},\dots,A_{n-s+1}\in\mathcal{K}_{\eta}^{*}(\mathbb{R})$ generating an hypercyclic abelian semigroup.

We repeatedly use the following multidimensional version of Kronecker's Theorem:
\
\\

{\bf Kronecker's Theorem} (See for example 
\cite{sh},\ Lemma 2.2, \cite{hAaM}, Lemma 4.2). Let $\alpha_{1},\dots, \alpha_{n}$ be
negative real numbers such that ~$1, \alpha_{1},\dots,
\alpha_{n}$ are linearly independent over $\mathbb{Q}$. Then the
set $$\mathbb{N}^{n} + \mathbb{N}[\alpha_{1},\dots,\alpha_{n}]^{T}
:=\left\{[s_{1},\dots,s_{n}]^{T} +
k[\alpha_{1},\dots,\alpha_{n}]^{T} :\  k, \;s_{1},\dots, s_{n}\in
\mathbb{N}\right\}$$ is dense in $\mathbb{R}^{n}$.

\begin{prop}\label{p:10} Let $n\in \mathbb{N}$, $n\geq 1$ and $s=0,1,\dots,n$. Then there exist ~$(n-s+1)$ vectors 
~$u_{1},\dots, u_{n-s+1}$ of ~$\mathbb{R}^{n}$ such that 
$$H_{s}: =\begin{cases}                                
\underset{k=1}{\overset{n-s+1}{\sum}}\mathbb{N}u_{k}+\underset{l=1}{\overset{s}{\sum}}
2\pi\mathbb{Z}f_{\eta}^{(l)} \ & \ \mathrm{ if } \ s\neq 0\\\\
\underset{k=1}{\overset{n+1}{\sum}}\mathbb{N}u_{k} \ & \ \mathrm{ if } \ s=0\end{cases}$$ 
is dense in $\mathbb{R}^{n}$.
\end{prop}

\begin{proof} \textit{Case $s\neq 0$}. Let $\alpha_{1},\dots, \alpha_{n}$ be negative real numbers such that ~$1, \alpha_{1},\dots, \alpha_{n}$ are linearly
independent over $\mathbb{Q}$. Define

$\mathcal{B}_{0}\backslash(e_{t_{1}},\dots, e_{t_{s}}):=
(e_{i_{s+1}},\dots, e_{i_{n}})$, where $e_{t_{l}} = f_{\eta}^{(l)}$, $l=1,\dots, s$ (see page 3) and define 
the matrix $S$ by
$$Se_{k} = \begin{cases}
                  2\pi f_{\eta}^{(k)} & \  \textrm{if} \ k=1,\dots, s, \\
                  e_{i_{k}} & \  \textrm{if} \ k= s+1,\dots, n.
                 \end{cases}$$
We see that $S \in \textrm{GL}(n;\mathbb{R})$. Set $u = [\alpha_{1},\dots, \alpha_{n}]^{T}$ and define 
\
\\
$$u_{k}:=\begin{cases}
                  Se_{s+k} & \ \textrm{if} \ \ k=1,\dots, n-s, \\
                  Su & \ \textrm{if} \ \ k = n-s+1.\\
                \end{cases}$$
Set $H^{\prime}_{s}:=\underset{k=1}{\overset{n-s}{\sum}}\mathbb{N}e_{s+k}+
\mathbb{N}u+\underset{l=1}{\overset{s}{\sum}}\mathbb{Z}e_{l}.$ We then have that
\begin{align*}
S(H^{\prime}_{s}) & =
\underset{k=1}{\overset{n-s}{\sum}}\mathbb{N}Se_{s+k}+
\mathbb{N}Su+\underset{l=1}{\overset{s}{\sum}}\mathbb{Z}Se_{l}\\ \
& =\underset{k=1}{\overset{n-s}{\sum}}\mathbb{N}u_{k}+
\mathbb{N}u_{n-s+1}+\underset{l=1}{\overset{s}{\sum}}2\pi\mathbb{Z}f_{\eta}^{(l)}\\
\ &
=\underset{k=1}{\overset{n-s+1}{\sum}}\mathbb{N}u_{k}+\underset{l=1}{\overset{s}{\sum}}2\pi\mathbb{Z}f_{\eta}^{(l)}\\
\ & = H_{s}
\end{align*}
Since $\mathbb{N}^{n}+\mathbb{N}u\subset H^{\prime}_{s}$, we see that $H^{\prime}_{s}$ is dense in
$\mathbb{R}^{n}$ by Kronecker's theorem, and thus so is $H_{s}$. 

\
\\
\textit{Case $s = 0$}. We let $$u_{k}= \begin{cases}
                  e_{k} & \ \textrm{if} \ \ 1\leq k \leq n, \\
                  u & \ \textrm{if} \ \ k = n+1.\\
                \end{cases}$$
\
\\
Then \;$H_{0}:=\underset{k=1}{\overset{n+1}{\sum}}\mathbb{N}u_{k} = \underset{k=1}{\overset{n}{\sum}}\mathbb{N}e_{k}+
\mathbb{N}u= \mathbb{N}^{n}+\mathbb{N}u$. So by Kronecker's theorem, $H_{0}$ is dense in
$\mathbb{R}^{n}$. 
This proves the proposition.
\end{proof}
\smallskip

\begin{proof}[Proof of Theorem ~\ref{T:4}] By Proposition ~\ref{p:10}, there exist $u_{1},\dots,u_{n-s+1}\in \mathbb{R}^{n}$
such that $H_{s}$ is dense in $\mathbb{R}^{n}$. 

\
\\
\textit{Case $s\neq 0$}. \\
 $-$ \textit{ Subcase: $r\neq 0$}. Set $u_{j} =
[u_{j,1},\dots, u_{j,r}; \widetilde{u}_{j,1},\dots,
\widetilde{u}_{j,s}]^{T}$ 

\
\\
with
$u_{j,k} = [x^{(j)}_{k,1},\dots,
x^{(j)}_{k,n_{k}}]^{T}$ and $\widetilde{u}_{j,l} =
[y^{(j)}_{l,1},-y^{\prime(j)}_{l,1}, \dots,
y^{(j)}_{l,m_{l}},-y^{\prime(j)}_{l,m_{l}}]^{T}, \\ j=1, \dots,
n-s+1, \ k=1, \dots, r, \ l=1,\dots,s$. 

\
\\
Let $B_{1}, \dots, B_{n-s+1}$ be defined by  $B_{j} = \textrm{diag}(B_{j,1},\dots,B_{j,r}; \
\widetilde{B}_{j,1},\dots,\widetilde{B}_{j,s})$, where

$$B_{j,k}=\left[\begin{array}{ccccc}
  x^{(j)}_{k,1} & \ & \ &\ & 0 \\
  \vdots & \ddots & \ & \  & \ \\
  \vdots & 0 & \ddots & \ & \ \\
  \vdots & \vdots & \ddots &\ddots & \ \\
  x^{(j)}_{k,n_{k}} & 0 & \dots & 0& x^{(j)}_{k,1}
\end{array}\right]\ \ \ \ \textrm{and}\ \ \ \ \widetilde{B}_{j,l}  = \left[\begin{array}{ccccc}
  C^{(j)}_{l,1} & \ & \ &\ & 0 \\
  \vdots & \ddots & \ & \  & \ \\
  \vdots & 0 & \ddots & \ & \ \\
  \vdots & \vdots & \ddots &\ddots & \ \\
  C^{(j)}_{l, m_{l}} & 0 & \dots & 0& C^{(j)}_{l,1}
\end{array}\right],$$
\
\\
with 
\
\\
\;$C^{(j)}_{l,i}=\left[\begin{array}{cc}
                               y^{(j)}_{l,i} & y^{\prime(j)}_{l,i} \\
                               -y^{\prime(j)}_{l,i} & y^{(j)}_{l,i}
                             \end{array}
\right], 1\leq j\leq n-s+1,\ 1\leq k \leq r,\  1\leq l \leq s, \ 1\leq i \leq m_{l}$.

\
\\
 As $B_{j,k}e_{\eta, k}=u_{j,k}$, then $B_{j}u_{\eta} = u_{j}$. 
\medskip

Let $G$ be the sub-semigroup of
$\mathcal{K}_{\eta}^{*}(\mathbb{R})$ generated by
$A_{1},\dots,A_{n-s+1}$, where  
\
\\
$$A_{j}=\textrm{diag}(A_{j,1},\dots,A_{j,r};
\widetilde{A}_{j,1},\dots,\widetilde{A}_{j,s}), \ j=1,\dots, n-s+1$$

$$A_{j,k} = \begin{cases}

                     e^{\frac{1}{2}B_{j,k}} \ &  \textrm{if} \ 1\leq k\neq j\leq r, \\
                    -e^{\frac{1}{2}B_{j,j}} \ &  \textrm{if} \  k=j.
           \end{cases}$$
and
 $$\widetilde{A_{j,l}} = e^{\frac{1}{2}\widetilde{B_{j,l}}}, \ \ l=1,\dots,s.$$

We have:
\
\\
\textbullet\; $r\leq n-s+1$ since $r+2s\leq n$,
\
\\
\textbullet\; $A^{2}_{j}=e^{B_{j}}$, $j=1,\dots,n-s+1$,
\
\\
\textbullet\; $A_{j,j}$ has a negative eigenvalue: $-e^{\frac{1}{2}x^{(j)}_{j,1}}$, for
every  $j=1,\dots, r$.
\
\\
As a consequence we have $\textrm{ind}(G)=r$.
\
\\
Firstly, we check that $G$ is abelian: For this, it suffices to show that
$A_{i}A_{j}=A_{j}A_{i}$ for every $1\leq i,j\leq n-s+1$, which is
equivalent to show that
 $B_{j}B_{j^{\prime}} = B_{j^{\prime}}B_{j}$ and 
$\widetilde{B}_{j}\widetilde{B}_{j^{\prime}} = \widetilde{B}_{j^{\prime}}\widetilde{B}_{j}$
for every $j, \ j^{\prime} = 1,\dots, n -s + 1$.

\
\\
Set $B_{j,k} = N_{j,k} + x^{(j)}_{k,1}I_{n_{k}}, \
\widetilde{B}_{j,l} = \widetilde{N}_{j,l} +
\widetilde{D}^{(j)}_{l,1},$ where 
 \
 \\
 $N_{j,k}= \left[\begin{array}{cc}
                  0 & 0 \\
                  T_{j,k} & 0
                \end{array}
\right]\in\mathbb{T}_{n_{k}}(\mathbb{R}) \ \textrm{with} \ 
T_{j,k} = \left[x^{(j)}_{k,2},\dots, x^{(j)}_{k,n_{k}}
\right]^{T}, \ k=1,\dots, r$
and \
\\
$\
\widetilde{D}^{(j)}_{l,1} = \textrm{diag}(C^{(j)}_{l,1},\dots,C^{(j)}_{l,1}), \ 
\widetilde{N}_{j,l}=\left[\begin{array}{cc}
                  0 & 0 \\
                  \widetilde{T}_{j,l} & 0
                \end{array}
\right]\in\mathbb{B}_{m_{l}}(\mathbb{R})$  with \ 
\
\\
$\widetilde{T}_{j,l} = \left[C^{(j)}_{l,2},\dots, C^{(j)}_{l,m_{l}}
\right]^{T}, \ l=1,\dots, s.$ We see that
$N_{j,k}N_{j^{\prime},k} = N_{j^{\prime},k}N_{j,k}=0$, for every
$j, j^{\prime}, k=1,\dots,s$. Hence $B_{j,k}B_{j^{\prime},k} =
B_{j^{\prime},k}B_{j,k}$. In the same way,
$\widetilde{N}_{j,k}\widetilde{N}_{j,k^{\prime}} =
\widetilde{N}_{j,k^{\prime}}\widetilde{N}_{j,k}=0$,
 $\widetilde{N}_{j,k}\widetilde{D}_{j,k^{\prime}}=\widetilde{D}_{j,k^{\prime}}\widetilde{N}_{j,k}$ and
 $\widetilde{N}_{j^{\prime},k}\widetilde{D}_{j,k}=\widetilde{D}_{j,k}\widetilde{N}_{j^{\prime},k}$
 for every $k=1,\dots,s$. 
Hence,
$\widetilde{B}_{j,k}\widetilde{B}_{j^{\prime},k} =
\widetilde{B}_{j^{\prime},k}\widetilde{B}_{j,k}$. We conclude that
$B_{j}B_{j^{\prime}} = B_{j^{\prime}}B_{j}$ and $\widetilde{B}_{j}\widetilde{B}_{j^{\prime}} = \widetilde{B}_{j^{\prime}}\widetilde{B}_{j}$.
\
\\
Secondly, by Proposition ~\ref{p:8} applying to $G^{2}$, we have that 
\begin{align*}
  \mathrm{g}_{\eta}^{2}(u_{\eta}) & = \underset{k=1}{\overset{n-s+1}{\sum}}\mathbb{N}B_{k}u_{\eta}+
\underset{l=1}{\overset{s}{\sum}}2\pi\mathbb{Z}f_{\eta}^{(l)} \\
  \ & = \underset{k=1}{\overset{n-s+1}{\sum}}\mathbb{N}u_{k}+
\underset{l=1}{\overset{s}{\sum}}2\pi\mathbb{Z}f_{\eta}^{(l)}\\ \ & = H_{s}
\end{align*}
By Proposition \ref{p:10}, it follows that 
$\overline{\textrm{g}_{\eta}^{2}(u_{\eta})} = \mathbb{R}^{n}$. Since $\textrm{ind}(G)=r$, so by Theorem ~\ref{T:2}, \ 
$\overline{G(u_{\eta})} = \mathbb{R}^{n}$.
\bigskip

$-$ \textit{Subcase $r=0$}. We consider the same construction above, where $G$ is the sub-semigroup of
$\mathcal{K}_{\eta}^{*}(\mathbb{R})$ generated by
$A_{1},\dots,A_{n-s+1}$, with  
\
\\
$$A_{j}=\textrm{diag}(\widetilde{A}_{j,1},\dots,\widetilde{A}_{j,s}), \ j=1,\dots, n-s+1$$ and 
$$\widetilde{A_{j,l}}=e^{\frac{1}{2}\widetilde{B_{j,l}}}, \ \ l=1,\dots,s.$$
In this case we have $\textrm{ind}(G)=0$. By running the same proof above, we get that $\overline{G(u_{\eta})} = \mathbb{R}^{n}$.
\medskip

\textit{Case $s= 0$}. Set $u = [\alpha_{1},\dots, \alpha_{n}]^{T}=
[v_{1},\dots, v_{r}]^{T}$ with 
$v_{k} = [x_{k,1},\dots,
x_{k,n_{k}}]^{T}$, and set $e_{j} =
[e_{j,1},\dots, e_{j,r}]^{T}$ 
with
$e_{j,k} = [e^{(j)}_{k,1},\dots,
e^{(j)}_{k,n_{k}}]^{T}$, where $e^{(j)}_{k,l} = \delta_{j}^{\sum\limits_{i=1}^{k-1}n_i +l}$($\delta_{p}^{q}$
is the Kronecker symbol),  
%$\begin{cases}
%                1 & \ \textrm{if} \ j= \sum\limits_{i=1}^{k-1}n_i +l \\
%                0 & \ \textrm{if} \ \  j\neq \sum\limits_{i=1}^{k-1}n_i +l  \end{cases}$ 
                $ l=1, \dots, n_k$,  $k=1, \dots, r$, $j=1, \dots,
n$. %In particular,  $e_{j,k} = e_{\eta,k}$ if $j= \sum\limits_{i=1}^{k-1}n_i +1$.
Let $B_{1}, \dots, B_{n+1}$ be defined by $B_{j} = \textrm{diag}(B_{j,1},\dots,B_{j,r})$, for every $j=1, \dots,
n+1$, where
\
\\

%for  $j\neq n+1$,  $B_{j,k} = \begin{cases}
%                  0_{n_{k}} & \ \textrm{if} \ \ k\neq j \\
%                  I_{n_{j}} & \ \textrm{if} \ \ k = j  \end{cases}$
%                  
for  $j\neq n+1$, $$B_{j,k} = \left[\begin{array}{ccccc}
 e^{(j)}_{k,1} & \ & \ &\ & 0 \\
 \vdots & \ddots & \ & \  & \ \\
\vdots & 0 & \ddots & \ & \ \\
\vdots & \vdots & \ddots &\ddots & \ \\
 e^{(j)}_{k,n_{k}} & 0 & \dots & 0 & e^{(j)}_{k,1}
\end{array}\right]$$

 In particular,  $$B_{j,k}= 
\begin{cases}  I_{n_{k}} & \ \textrm{if} \ \ j= \sum\limits_{i=1}^{k-1}n_i +1\\
                  0_{n_{k}} & \ \textrm{if} \ \ j\notin \{\sum\limits_{i=1}^{k-1}n_i +1, \dots, \sum\limits_{i=1}^{k}n_i\}
                  \end{cases}$$
\
\\
and for  $j= n+1$, $$B_{j,k} = \left[\begin{array}{ccccc}
  x_{k,1} & \ & \ &\ & 0 \\
  \vdots & \ddots & \ & \  & \ \\
  \vdots & 0 & \ddots & \ & \ \\
  \vdots & \vdots & \ddots &\ddots & \ \\
  x_{k,n_{k}} & 0 & \dots & 0& x_{k,1}
\end{array}\right]$$

\
\\
As $$B_{j,k}e_{\eta, k} = 
\begin{cases}  e_{j,k} & \ \textrm{if} \ \ j\neq n+1 \\
                  v_{k} & \ \textrm{if} \ \  j=n+1  \end{cases}$$
\
\\
we obtain $$B_{j} u_{\eta} = \begin{cases} e_{j} & \ \textrm{if} \ \ j\neq n+1 \\
                  u & \ \textrm{if} \ \  j=n+1
 \end{cases}$$

\
\\
Let $G$ be the sub-semigroup of
$\mathcal{K}_{\eta}^{*}(\mathbb{R})$ generated by
$A_{1},\dots,A_{n+1}$, where $\eta$ is of length $r$ and 
$A_{j}=\textrm{diag}(A_{j,1},\dots,A_{j,r}), \ j=1,\dots, n+1$,
with $$A_{j,k}  = \begin{cases} e^{\frac{1}{2}B_{j,k}} \ &  \ \textrm{if} \ 1\leq k\neq j\leq r, \\
                  -e^{\frac{1}{2}B_{j,j}} \ & \ \textrm{if} \ \ k= j   
 \end{cases}$$

\
\\
We have:

\
\\
\textbullet\; $A^{2}_{j} = e^{B_{j}}, \ j=1,\dots,n+1$,
\
\\
\textbullet\; $A_{j,j}$ has a negative eigenvalue: $-e^{\frac{1}{2}e^{(j)}_{j,1}}$, \;for
every  $j=1,\dots, r$.
\
\\
 Then by above, we have $\textrm{ind}(G)=r$. Moreover it is easy to check that $G$ is abelian. 
 Now apply Proposition ~\ref{p:8} to $G^{2}$, we have that 
\begin{align*}
  \mathrm{g}^{2}_{u_{\eta}}  = \underset{j=1}{\overset{n+1}{\sum}}\mathbb{N}B_{j}u_{\eta}
  = \underset{j=1}{\overset{n}{\sum}}\mathbb{N}e_{j}+ \mathbb{N}u = H_{0}.
\end{align*}
\
\\
By Proposition \ref{p:10}, it follows that 
$\overline{\textrm{g}_{\eta}^{2}(u_{\eta})} = \mathbb{R}^{n}$. Since $\textrm{ind}(G)=r$, so by Theorem ~\ref{T:2}, \ 
$\overline{G(u_{\eta})} = \mathbb{R}^{n}$.
\end{proof}
\medskip

{\it Proof of Corollary \ref{C:6n}}. Take
$s=n-\left[\frac{n+1}{2}\right]$. Then $s\in\mathbb{N}$, and by Theorem \ref{T:4}, there exist
$(\left[\frac{n+1}{2}\right]+1)$ matrices in
$\mathcal{K}_{\eta}^{*}(\mathbb{R})$ where $\eta$ has length $r+2s$, generating a
hypercyclic abelian semigroup. We conclude, by Corollary \ref{C:4}, that
$(\left[\frac{n+1}{2}\right]+1)$ is the minimal number of matrices
having such property. \qed
\medskip

\section{\bf Example}
\medskip

\begin{exe} Let $G$ be the semigroup generated by $A_{1} = \mathrm{diag}(e^{\pi}, e^{\pi})$,
$A_{2} =\left[\begin{array}{cc}
            -1 & 0 \\
            -\pi & -1
          \end{array}
\right]$ and $A_{3} =e^{-\pi\sqrt{2}}\left[\begin{array}{cc}
            1 & 0 \\
            -\pi\sqrt{3} & 1
          \end{array}
\right].$

\
\\
Then $G$ is abelian and hypercyclic.
\end{exe}

\begin{proof} By construction, $G$ is an abelian sub-semigroup of $\mathbb{T}^{*}_{2} (\mathbb{R})$ with $\mathrm{ind}(G)=1$,
$\eta = (2)$, $u_{\eta} = e_{1}$ and $A^{2}_{k} = e^{B_{k}}$, $k = 1,2,3$, where
$B_{1} = \mathrm{diag}(2\pi; 2\pi)$, $B_{2}
=\left[\begin{array}{cc}
            0 & 0 \\
            2\pi & 0
          \end{array}
\right]$ and $B_{3} =\left[\begin{array}{cc}
            -2\pi\sqrt{2} & 0 \\
            -2\pi \sqrt{3}& -2\pi\sqrt{2}
          \end{array}
\right].$
\
\\
By Proposition ~\ref{p:8}, we have that
%\begin{align*}
$\mathrm{g}_{\eta}^{2}(e_{1}) =\underset{k=1}{\overset{3}{\sum}}\mathbb{N}B_{k}e_{1} =
2\pi H,$
%\end{align*}
where 
\
\\
$H := \mathbb{N}e_{1} + \mathbb{N}e_{2}  +
\mathbb{N}[-\sqrt{2}, -\sqrt{3}]^{T}= \mathbb{N}^{2}  +
\mathbb{N}[-\sqrt{2}, -\sqrt{3}]^{T}$. By Kronecker's Theorem,
$\overline{H}=\mathbb{R}^{2}$ since $1, -\sqrt{2}$ and $-\sqrt{3}$
are linearly independent over $\mathbb{Q}$. Therefore, by Theorem ~\ref{T:2}, $\overline{G(e_{1})}=
 \mathbb{R}^{2}$.
\end{proof}

\bibliographystyle{amsplain}
\vskip 0,4 cm
%%% ----------------------------------------------------------------------

\end{document}